\documentclass[ sigconf,nonacm,table]{acmart}

\copyrightyear{2025}
\acmYear{2025}
\setcopyright{cc}
\setcctype{by}
\acmConference[KDD '25]{Proceedings of the 31st ACM SIGKDD Conference on Knowledge Discovery and Data Mining V.2}{August 3--7, 2025}{Toronto, ON, Canada}
\acmBooktitle{Proceedings of the 31st ACM SIGKDD Conference on Knowledge Discovery and Data Mining V.2 (KDD '25), August 3--7, 2025, Toronto, ON, Canada}
\acmDOI{10.1145/3711896.3737055}
\acmISBN{979-8-4007-1454-2/2025/08}

\usepackage{float}
\usepackage{booktabs}       
\usepackage{scalefnt}
\usepackage{nicefrac}       
\usepackage{microtype}      
\usepackage{optidef}
\usepackage{amsmath,amsfonts}
\usepackage{bbm}
\usepackage{mathtools}
\usepackage{stmaryrd}
\usepackage{subcaption}
\usepackage{graphicx} 
\usepackage{bbm}
\usepackage{paralist}
\usepackage{enumitem}
\setitemize{noitemsep,topsep=0pt,parsep=0pt,partopsep=0pt}

\newtheorem{proposition}{Proposition}
\newtheorem{corollary}{Corollary} 

\usepackage{svg}
\usepackage{tabularray}
\usepackage{supertabular}
\usepackage{booktabs} 

 \usepackage[linesnumbered,ruled]{algorithm2e}

\SetKwInput{KwInput}{Input}                
\SetKwInput{KwOutput}{Output}  
\SetKwRepeat{Repeat}{repeat}{until}
\SetKwRepeat{RepeatFor}{repeat}{for}
\SetKwRepeat{Repeat}{repeat}{until}

\SetAlCapNameFnt{\footnotesize}
\SetAlCapFnt{\footnotesize}

\title{MOSS: Multi-Objective Optimization for Stable Rule Sets}

\author{Brian Liu}
\affiliation{%
  \institution{Massachusetts Institute of Technology }
   \city{Cambridge}
   \state{Massachusetts}
  \country{USA}}
\email{briliu@mit.edu}

\author{Rahul Mazumder}
\affiliation{%
  \institution{Massachusetts Institute of Technology }
   \city{Cambridge}
   \state{Massachusetts}
  \country{USA}}
\email{rahulmaz@mit.edu}
 
\begin{document}

\begin{abstract}

We present MOSS, a multi-objective optimization framework for constructing stable sets of decision rules. MOSS incorporates three important criteria for interpretability: sparsity, accuracy, and stability, into a single multi-objective optimization framework. Importantly, MOSS allows a practitioner to rapidly evaluate the trade-off between accuracy and stability in sparse rule sets in order to select an appropriate model. We develop a specialized cutting plane algorithm in our framework to rapidly compute the Pareto frontier between these two objectives, and our algorithm scales to problem instances beyond the capabilities of commercial optimization solvers. Our experiments show that MOSS outperforms state-of-the-art rule ensembles in terms of both predictive performance and stability.
\end{abstract}



\maketitle

\section{Introduction}

Rule sets are prized in high-stakes applications of machine learning, such as criminal justice and healthcare operations, for balancing a high degree of transparency with good predictive performance \cite{rudin2014algorithms,han2014rule}. These ensembles consist of decision rules, or sequences of if-then antecedents (i.e. splits) that partition a datasets and assign a prediction to the partition; summing these predictions produces the output of the ensemble \cite{friedman2008predictive}. Sparse rule ensembles, of a manageable set of less than 20 or so decision rules \cite{liu2023fire}, are especially useful in high-stakes applications since these rules can be audited by hand for fairness or bias concerns. However, sparse rule sets should also be stable and accurate to be considered trustworthy \cite{yu2018three}. Intuitively, an algorithm that produces vastly different rule sets across small data perturbations is unlikely to be trusted. Likewise, rule sets that have insufficient explanatory power, or that perform poorly out-of-sample, should not be relied upon\footnote{We direct readers to the Predictability, Computability, and Stability framework for veridical data science\cite{yu2018three,murdoch2019definitions,yu2020veridical} for an in-depth discussion on criteria for trustworthiness in interpretable machine learning.}.

Sparsity, accuracy, and stability are three competing objectives to consider when constructing interpretable  models. The trade-off between model sparsity and accuracy has been well-explored, and popular frameworks such as \textsc{glmnet} and \textsc{L0Learn} allow practitioners to rapidly evaluate how model size impacts predictive performance \cite{hazimeh2023l0learn,hastie2021introduction}. The trade-off between model stability and accuracy, on the other hand, is poorly understood. Frameworks that improve model stability, such as stability selection \cite{meinshausen2010stability}, do not account for predictive accuracy when constructing models. We illustrate this in the motivating example below, where we use a stability selection-based framework to construct stable rule sets \cite{meinshausen2010stability,benard2021sirus}.

The general stability selection framework for constructing rule sets proceeds as follows. Given a dataset, we apply a  sampling procedure such as the bootstrap to generate multiple perturbed datasets. We construct a set of decision rules on each perturbation; one such algorithm, SIRUS, does so by fitting a single decision tree and taking the leaf nodes as decision rules \cite{benard2021sirus}. Across all perturbations, we compute the frequency at which unique decision rules appear, and we return the set of rules with frequencies above some threshold. By extracting rules that are consistently constructed across data perturbations, stability selection improves model stability. However, the framework does not select rules with respect to their accuracy on the original data. As such, the predictive performance of the rule sets may suffer.

\begin{figure}
    \centering
    \includegraphics[width = 0.475\textwidth]{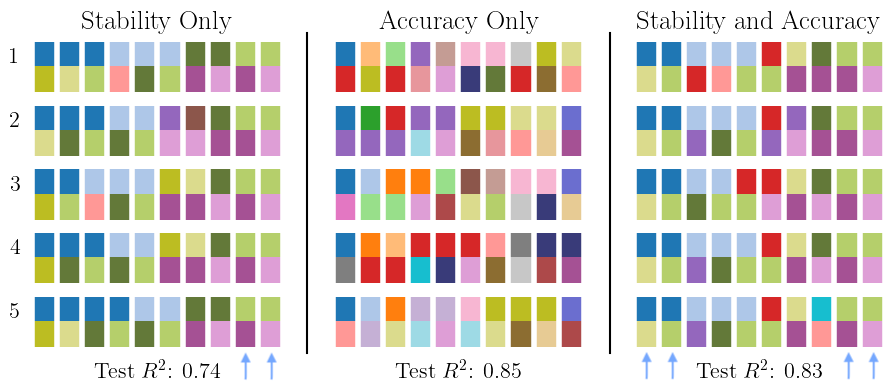}
    \caption{Rule sets constructed across 5 train-test splits. The blue arrows show rules that are stable across all splits.}

    \label{motivating_example.fig}
\end{figure}

In the left panel in Figure \ref{motivating_example.fig}, we apply stability selection (SIRUS) to construct 9 decision rules of interaction depth 2, across 5 different train-test partitions of the \textsc{Galaxy} regression dataset \cite{OpenML2013}. Each rectangle in the panel represents a decision rule and contains two vertically-stacked squares to represent each split. The squares are color-coded with respect to split locations, i.e, squares with the same color are identical splits with respect to features and values. We observe that two rules are consistently recovered across partitions (indicated by the blue arrows) and that the average out-of-sample $R^2$ of the constructed rule set is \textbf{0.74}. In the center panel in Figure \ref{motivating_example.fig}, we use  RuleFit \cite{friedman2008predictive} to construct 9 decision rules with respect to only predictive accuracy. We observe that the constructed rules perform much better, with an average out-of-sample $R^2$ of \textbf{0.85}, however, the rules are unstable. No rules are consistently constructed across all train-test splits. From this example, we see that there is a trade-off between stability and accuracy that should be explored when constructing sparse rule sets.

In this paper we introduce MOSS, a multi-objective optimization framework for constructing stable rule sets. In contrast to existing algorithms, such as SIRUS and RuleFit, MOSS incorporates stability, sparsity, and accuracy into a single unified optimization problem. Importantly, we develop a novel optimization formulation and algorithm in MOSS that allows our framework to efficiently explore the Pareto frontier between accuracy and stability. As such, practitioners can use MOSS to  rapidly evaluate the trade-off between these two objectives in order to select an appropriate model. Our experiments show that MOSS can extract sparse stable models that outperform state-of-the-art  rule-based algorithms. In fact, when we apply MOSS to our motivating example discussed above we construct rule sets more stable than SIRUS with better out-of-sample predictive performance. In the right panel in Figure \ref{motivating_example.fig} we show the rule sets constructed by MOSS. We observe that 4 rules are consistently constructed across all train-test splits and that the rule sets have an average out-of-sample $R^2$ of \textbf{0.83}. We summarize the contributions of our paper below.

\begin{itemize}[leftmargin=5mm]
    \item We develop a multi-objective optimization framework for constructing rule sets that jointly accounts for accuracy, stability, and sparsity (\S\ref{moss1.section}).
    \item We develop specialized optimization algorithms (\S\ref{optalg.sec}) that allow practitioners to efficiently explore the Pareto frontier between model accuracy and stability (\S\ref{efficient_pareto.section}).
    \item We show that MOSS can outperform state-of-the-art rule-based algorithms in terms of both stability and accuracy (\S\ref{experiments.section}).
\end{itemize}

The remainder of our paper is organized as follows. We first discuss some preliminaries on model stability and the stability selection framework. We show that stability selection can be re-interpreted from an optimization viewpoint, which motivates our exploration into multi-objective optimization for both stability and accuracy. We then present the MOSS framework, along with the specialized optimization algorithm we use to solve problems in MOSS. We conclude with our experimental results, where we evaluate the stability and predictive performance of rule sets created using MOSS against various state-of-the-art algorithms, followed by a discussion connecting MOSS to related works.

\section{Problem Formulation}

We first overview why stability is crucial in trustworthy and interpretable machine learning and introduce how we assess the stability of rule sets. Following this, we discuss stability selection, a popular statistical framework for improving model stability. We show how an optimization reinterpretation of this framework motivates our formulation of MOSS.


\subsection{Model Stability} \label{model_stability_def.section}
Stability has long been considered to be a prerequisite for trustworthiness when interpreting models \cite{yu2013stability,yu2020veridical}. At the minimum, reliable conclusions drawn from models must be replicable across repeated analyses, a foundational idea in scientific discovery \cite{popper2005logic}. From a statistical perspective, we are interested in the stability of the results of repeated analyses conducted across datasets that differ by reasonable perturbations \cite{yu2013stability}. We formalize this in the context of constructing rule sets below.

We want to assess if an algorithm constructs the same set of rules when applied across different data perturbations. Say that we have datasets $X_i$ and $X_j$ that differ by some reasonable perturbation, for example, these two datasets may come from different training-test partitions of the original dataset $X$. Using the same algorithm, we construct rule sets $R_i$ and $R_j$ on each dataset. To measure the similarity between the two rule sets we use the Dice-Sorensens coefficient, defined by:
\begin{equation}
   \text{DSC}(R_i,R_j) = \frac{2|R_i \cap R_j|}{|R_i|+|R_j|}.
\end{equation}
This coefficient captures the proportion of decision rules shared between rule sets and is commonly used in assessing the stability of statistical models \cite{nogueira2018stability}.

To assess the stability of our rule algorithm, we generate many perturbed datasets $X_1 \ldots X_T$ and apply our algorithm to construct rule sets $R_1 \ldots R_T$. We report the average pairwise Dice-Sorensens coefficient across all rule sets as the \emph{empirical stability} of our algorithm, given by:
\begin{equation}
    \text{Empirical Stability} = \frac{T(T-1)}{2}\sum_{i\neq j} \text{DSC}(R_i,R_j).
\end{equation}
This metric captures the average proportion of rules shared between rule sets constructed across data perturbations. As an aside, in  \S\ref{choice_of_stability.section} of the appendix, we discuss alternative metrics for rule stability and demonstrate that the results of our paper remain consistent across these metrics.

In the next section, we discuss how stability selection framework can be used to improve the stability of rule-based algorithms.

\subsection{Stability Selection}
The goal of stability selection is to remove the unstable components of a model. Given training dataset $X \in \mathbb{R}^{n\times p}$, the procedure starts by generating perturbed datasets $X_1 \ldots X_B$ by sub-sampling rows or the bootstrap. On each dataset, we apply the sample algorithm to generate rule sets $R_1 \ldots R_B$. Let $m$ denote the total number of \emph{unique} rules constructed across all sets. For each unique rule $r_i$ for $i \in [m]$ we compute the proportion of rule sets in which $r_i$ appears. We store these selection proportions in vector $\Pi \in R^m$. Stability selection extracts the rules with selection proportions above threshold $\tau$, $\{r_i | \Pi_i > \tau\}$; these rules are considered the stable components of the model. By pruning the unstable components, stability selection improves the empirical stability of the final model.

As an aside, the SIRUS algorithm mentioned in the introduction cleverly applies stability selection to the construction of rule sets by leveraging random forests \cite{benard2021sirus}. Random forests consist of decision trees fit on bootstrapped datasets, and the leaves of each decision tree form a set of decision rules. SIRUS extracts the rules that appear consistently across trees in the random forest.

 Threshold $\tau$ is a hyperparameter that controls the sparsity $k$ of the selected model. For rule sets, $k$ is often prespecified to be a managable number of rules (at most 15-20) to preserve interpretability \cite{liu2023fire,benard2021sirus,wei2019generalized}. Next, we show that we can reinterpret stability selection as an optimization problem for a fixed model size $k$.

\subsection{Optimization Reinterpretation}\label{knapsack.section}

Stability selection eliminates unstable rules by discarding those with selection proportions below threshold $\tau$. For a predetermined model size of $k$ rules, discarding unstable rules can be achieved by selecting rules with the top-$k$ selection proportions $\Pi$. We formalize this task below.

Given a set of $r_1 \ldots r_m$ unique decision rules with selection proportions $\Pi \in \mathbb{R}^m$, we introduce binary decision variables $z_i \in \{0,1\}$, for $i \in [m]$, to indicate if rule $r_i$ is selected. Selecting the top-$k$ rules with the largest selection proportions can be expressed by this optimization problem:
\begin{align}\label{knapsack.problem}
\begin{split}
    \max_{z_1, \ldots, z_m} H_1(z) = \sum_{i=1}^m \Pi_i z_i
    \end{split}
    \begin{split}
        \text{s.t.} & \quad {\sum_{i=1}^m }  z_i \leq k, \quad z \in \{0,1\}^m.
    \end{split}
\end{align}
We denote objective $H_1(z)$ as the \emph{in-sample stability} of solution $z$, i.e., the sum of selection proportions of the selected rules. For a fixed model size $k$, stability selection maximizes in-sample stability as a proxy for the empirical stability of the model.

Problem \ref{knapsack.problem} is a binary knapsack optimization problem and optimal solution $z^*$ can be trivially computed by taking the $k$-largest elements of $\Pi$. Moreover, we can also compute a sequence of $t$ solutions of progressive suboptimality by sorting $\Pi$ in descending order sliding a window of length $k$ along the sorted vector. This yields a sequence of solutions $z^*, z^1, z^2, \ldots, z^t$ with decreasing in-sample stability objectives $H_1(z^*) > H_1(z^1) > H_1(z^2) \ldots > H_1(z^t)$. Below, we show that exploring these sequences of suboptimal solutions yield important insights towards the trade-off between model accuracy and stability, insights that motivate the development of our MOSS framework.

\subsection{Accuracy of Suboptimal Solutions} \label{subopt_accuracy.section}

In this section, we examine a sequence of solutions $z^*, z^1, \ldots, z^t$ to Problem \ref{knapsack.problem} that are progressively suboptimal with respect to in-sample stability. Our goal is to explore how these solutions perform in terms of model accuracy, which we define for solution $z$ below.

\textbf{Model Accuracy:} Given a set of decision rules $r_1 \ldots r_m$ where $f_i(X) \in \mathbb{R}^n$ is the prediction of rule $r_i$, solution vector $z \in \{0,1\}^m$ indicates which rules are selected.
We want to assess how well the predictions of this selected rule set fit the data. Predictions of a rule set are determined by a linear combination of the rules  $\sum_{i = 1}^m w_i f_i(X)$, where decision variable $w_i  \neq 0$ if and only if $z_i \neq 0$, for all $i \in [m]$. Given response $y \in \mathbb{R}^m$ and solution vector $z$, we use regularized in-sample training loss to measure model accuracy, which is given by: \begin{align}\label{in_sample_accuracy.problem}
H_2(z) = \min_{w_1, \ldots, w_n}  \frac{1}{2}\Vert y -  \sum_{i=1}^m f_i(X) w_i\Vert_2^2 + \frac{1}{2\gamma} \Vert w \Vert_2^2, \\
\text{s.t.} \quad w_i(1-z_i) = 0 \quad  \forall \ i \in [m], \quad z \in \{0,1 \}^m. \notag
\end{align}
We add the ridge penalty for computational reasons, however, we note that the regularizer may improve performance in noisy settings \citep{mazumder2023subset}. Hyperparameter $\gamma > 0$ controls the degree of shrinkage over $w$. Note that for any $z$, the optimal solution to the minimization problem in expression \ref{in_sample_accuracy.problem}, $w^*$, can be obtained in closed-form through a back-solve. We show later in \S\ref{BIP.section} that we can rewrite expression \ref{in_sample_accuracy.problem} in terms of just $z$.


Now that we have defined a way to assess model accuracy, through in-sample loss $H_2(z)$, we can begin our exploration. Using the procedure discussed in \S\ref{knapsack.section}, we apply Problem \ref{knapsack.problem} to construct rule sets on the \textsc{ESL} dataset \cite{OpenML2013}. We start with a collection of $10^3$ decision rules generated from a random forest and, for a fixed set of $k = 20$ decision rules, we compute the stability selection optimal solution $z^*$ and a sequence of 10 progressively suboptimal solutions $z^1, z^2, \ldots, z^{10}$. For each solution, we compute in-sample stability $H_1(z)$ and in-sample loss $H_2(z)$. 

\begin{figure}
    \centering
      \includegraphics[width=0.4\textwidth]{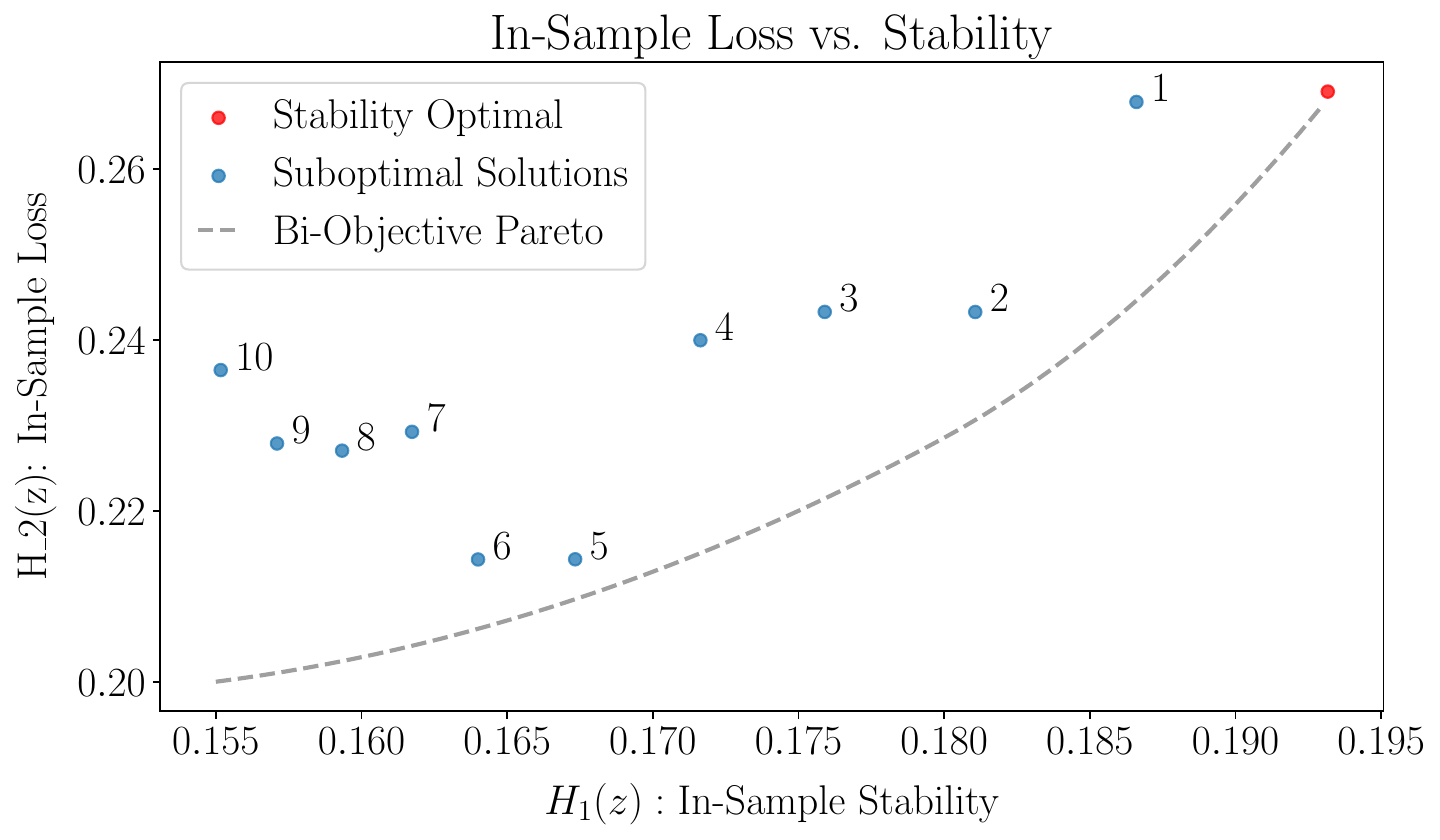}
    \caption{We apply stability selection to extract compact rule sets from $10^3$ candidate rules fit on the ESL dataset \cite{OpenML2013}. Solutions suboptimal in terms of stability have much lower in-sample loss. The goal of MOSS is to compute the bi-objective Pareto front between accuracy and stability (gray line).  }
\label{knapsack_explore.fig}
\end{figure}


In Figure \ref{knapsack_explore.fig}, we plot in-sample loss against in-sample stability for all of the solutions. The red point in the top right of the plot shows the optimal stability selection solution $z^*$ and the blue points show suboptimal solutions $z^1 \ldots z^{10}$. We observe from this plot that solutions that are slightly suboptimal with respect to in-sample stability may have significantly lower in-sample loss.

This observation motivates MOSS. The goal of our framework is to efficiently compute the Pareto frontier between in-sample stability $H_1(z)$ and in-sample loss $H_2(z)$, as hypothesized by the gray dashed line in Figure \ref{knapsack_explore.fig}. Importantly, we demonstrate later in \S\ref{experiments.section} that choosing solutions along this Pareto that are suboptimal with respect to in-sample stability can significantly improve the out-of-sample predictive performance of the model without compromising out-of-sample empirical stability of the model. We present our MOSS framework below.

\section{MOSS Framework} \label{MOSS_framework.section}
In this section we present the MOSS framework. We first discuss our multi-objective optimization formulation that accounts for both stability and accuracy. The optimization programs in this formulation are NP-hard, so we develop a specialized cutting plane-based optimization algorithm to solve problems in MOSS efficiently. Moreover, we develop a novel technique that exploits problem structure to efficiently compute the entire Pareto frontier between stability and accuracy.

\subsection{Multi-Objective Optimization Formulation} \label{moss1.section}

Given data $X \in \mathbb{R}^{n \times p}$ and target $y \in \mathbb{R}^n$, the goal of MOSS is to construct rule sets that are both accurate and stable. We first follow the stability selection (SIRUS) framework and construct a large collection of unique candidate rules $r_i$ for $i \in [m]$, with selection proportions $\Pi$. Our goal is to select which candidate rules to include into the final rule set, so we let $z \in \{0,1\}^m$ represent a vector of binary decision variables that indicate which rules are selected.

Recall that in \S\ref{knapsack.section} and \S\ref{subopt_accuracy.section}, we define in-sample stability $H_1(z)$ as a proxy for empirical stability and in-sample loss $H_2(z)$ as a proxy for out-of-sample accuracy. The goal of MOSS is to fit $z$ to \emph{maximize} in-sample stability and \emph{minimize} in-sample loss. We express maximizing $H_1(z)$ as minimizing $-H_1(z)$ and formalize this in the problem below. MOSS uses the following bi-objective integer program to select a rule set of $k$ decision rules:\begin{align}
\begin{split}
    \min_z \ -H_1(z), \quad  H_2(z) 
    \end{split}
    \begin{split}
        \text{s.t.} & \quad  {\sum_{i=1}^m }  z_i \leq k, \quad  z \in \{0,1\}^m.
    \end{split}
\end{align}

The goal of our framework is to compute the Pareto frontier between objectives $H_1(z)$ and $H_2(z)$. We apply the $\epsilon$-constraint method to accomplish this.

\subsubsection{ $\epsilon$-Constraint Method:} The $\epsilon$-constraint method computes the Pareto frontier of a bi-objective optimization problem by solving a sequence of single-objective problems. The method involves moving one objective into the constraint set and constraining that objective to be more extreme than a fixed value $\epsilon$. We sweep through values of $\epsilon$ while solving the corresponding single-objective problem to compute the Pareto frontier. This method is advantageous in that it can recover non-convex Pareto frontiers \cite{emmerich2018tutorial}.

For MOSS, we move stability objective $-H_1(z)$ into the constraint set and constrain $-H_1(z)$ to be less than or equal to $-\epsilon$, where $\epsilon$ is non-negative.  We sweep through values of $\epsilon$ while solving this single-objective integer program:
\begin{align}
\begin{split}
    \min_z \ H_2(z) 
    \end{split}
    \begin{split}
        \text{s.t.} &  \quad H_1(z) \geq \epsilon, \quad  {\sum_{i=1}^m }  z_i \leq k,  \quad z \in \{0,1\}^m.
    \end{split}
\end{align}
to compute our Pareto frontier. By plugging in our definition for model accuracy objective $H_2(z)$, expression \ref{in_sample_accuracy.problem}, we can express this problem as the following mixed integer program:


\begin{equation}\label{mainoptimizationproblem_temp}
\begin{alignedat}{2}
&\min_{w, z} &&\quad \frac{1}{2}\Vert y -  \sum_{i=1}^m f_i(X) w_i\Vert_2^2 + \frac{1}{2\gamma} \Vert w \Vert_2^2 , \\
&\text{s.t.} &&\quad \sum_{i=1}^m \Pi_i z_i \geq \epsilon, \quad \sum_{i=1}^m z_i \leq k, \\
&&& \quad w_i(1-z_i) = 0 \quad  \forall \ i \in [m],  \quad z \in \{0,1\}^m.
\end{alignedat}
\end{equation}

 We show below that we can reformulate this optimization problem entirely in terms of $z$.

\subsubsection{Binary Integer Reformulation:} \label{BIP.section} In Proposition 1, we show that we can reformulate Problem \ref{mainoptimizationproblem_temp} into a binary integer program by re-expressing the objective in terms of $z$.

\begin{proposition} Let prediction matrix $M \in \mathbb{R}^{n \times m}$ contain predictions $f_i(X) \in \mathbb{R}^n$ in column $i$, for $i \in [n]$. We can re-write our expression for regularized in-sample  loss, $H_2(z)$, given solution vector $z$ as:
 \begin{equation}
     H_2(z) = \frac{1}{2}y^{\intercal}\biggl(\mathbb{I}_n +\gamma \sum_{i=1}^m z_i M_i M_i^{\intercal}\biggr)^{-1} y,
 \end{equation}
 where $M_i$ is the $i$-th column of $M$.
 \end{proposition}

  We show this derivation in the appendix (\ref{prop1.appx}). As such, we can rewrite Problem \ref{mainoptimizationproblem_temp} as: \begin{equation}\label{mainoptimizationproblem}
\begin{alignedat}{2}
&\min_{ z} &&\quad \frac{1}{2}y^{\intercal}\biggl(\mathbb{I}_n +\gamma \sum_{i=1}^m z_i M_i M_i^{\intercal}\biggr)^{-1} y, \\
&\text{s.t.} &&\quad \sum_{i=1}^m \Pi_i z_i \geq \epsilon , \quad \sum_{i=1}^m z_i \leq k, \quad z \in \{0, 1\}^m.
\end{alignedat}
\end{equation}


This binary integer program is the main optimization problem that we use to construct rule sets in MOSS, where $\epsilon$ and $k$ are non-negative hyperparameters that control the stability and sparsity of the model respectively.

Problem \ref{mainoptimizationproblem} has $m$ decision variables and $m+2$ constraints; $m$ represents the size of the candidate rule set, which can be very large. For example, the SIRUS approach uses random forests to generate thousands of candidate rules. As such, instances of Problem \ref{mainoptimizationproblem} are often intractable and exceed the capabilities of off-the-shelf optimization software, as we show in \S\ref{computation_time_eval.section}. We can, however, exploit problem structure to develop a tailored optimization algorithm to solve this problem efficiently. Our specialized algorithm hinges on the fact that we move stability objective $H_1(z)$ into the constraint set, and not accuracy objective $H_2(z)$. By moving $H_1(z)$, we obtain constraints that are linear with respect to $z$. Also, the combinatorial space of feasible solutions for $z$ shrink as $\epsilon$ increases. We exploit these properties in our algorithm below.

\subsection{Cutting Plane Algorithm}\label{optalg.sec}

Here, we develop a specialized cutting plane algorithm to efficiently solve Problem \ref{mainoptimizationproblem} to optimality. Our algorithm leverages the fact that while integer programs with nonlinear objectives like Problem \ref{mainoptimizationproblem} are often intractable, integer linear programs (ILPs) with linear objectives can be efficiently solved using off-the-shelf methods \cite{bertsimas2020sparse}. Therefore, we reformulate the task of solving Problem \ref{mainoptimizationproblem} into solving a sequence of ILPs. Problem \ref{mainoptimizationproblem} is a binary integer program with linear constraints and we show in the proposition below that objective is convex.

\begin{proposition} Function $H_2(z)$ is convex on domain $z \in \mathbb{R}^m$ such that $z_i \in [0,1]$, for all $i \in [m]$.
\end{proposition}

We show the full proof of this proposition in the appendix (\ref{prop2.appx}). Since $H_2(z)$ is convex, we can apply an cutting plane algorithm in order to solve Problem \ref{mainoptimizationproblem} efficiently \cite{duran1986outer,bertsimas2020sparse}. Let $\nu$ represent the lower bound on objective value $H_2(z)$. At each iteration $t$ of our cutting plane algorithm, we first add cutting plane 
$\nu \geq H_2(z^t) + \nabla H_2(z^t)(z-z^t)$. The approximation function at iteration $t$ is given by $\max_{j \in [t]} H_2(z^j) + \nabla H_2(z^j)(z-z^j)$. We then minimize this function, with respect to the linear constraints in Problem \ref{mainoptimizationproblem} to update $\nu$ and $z$ for the subsequent iteration. This minimization involves solving an integer linear program, which can be done efficiently using off-the-shelf methods. We terminate our algorithm when $ H_2(z^t) \leq \nu$ and our algorithm terminates after a finite number of iterations and returns the optimal solution to Problem \ref{mainoptimizationproblem}. We discuss convergence further in \S\ref{alg_convergence.section} of the appendix. Our full cutting plane algorithm is presented below.

  \begin{algorithm}[h]
    \caption{ Cutting Plane}\label{OA.alg}
    \footnotesize
    $z^0 \leftarrow$ \text{warm start}

    $\nu^0 \leftarrow H_2(z^0), \ t \leftarrow 0$

    \While{$H_2(z^t) >  \nu $}{
    add constraint $\nu \geq H_2(z^t) + \nabla H_2(z^t)(z-z^t)$ to ILP (Problem \ref{ILP_problem})
    
    solve ILP (\ref{ILP_problem}) for $\nu^{t+1}, z^{t+1}$

    t = t+1
    }
    \textbf{return} $z^*$
  \end{algorithm}


\begin{mini!}|s|
{\nu, z} {\nu \label{ILP_obj}}{\label{ILP_problem}}{}
\addConstraint{ \sum_{i=1}^m \Pi_i z_i \geq \epsilon  \label{ILPc1}}
\addConstraint{ \sum_{i=1}^m  z_i\leq k, \quad z \in \{0, 1\}^m \label{ILPc2}}
\addConstraint{
\begin{aligned}
\nu &\geq H_2(z^j) + \nabla H_2(z^i)(z-z^j),
&\quad \forall \ j \in [t].\label{ILPc3}
\end{aligned}
}
\end{mini!}

\vspace{5mm}

In each iteration of our cutting plane algorithm (Algorithm \ref{OA.alg}), we must solve integer linear program Problem \ref{ILP_problem}. Off-the-shelf solvers can typically solve this problem in seconds, for problem sizes on the order of $m \sim 10^5$ decision variables. Below, we discuss two attributes of Algorithm \ref{OA.alg} that further improve computational efficiency.

\subsubsection{Sandwiching Constraints:} 
\begin{figure}[h]
  \begin{minipage}[c]{0.25\textwidth}
    \includegraphics[width=\textwidth]{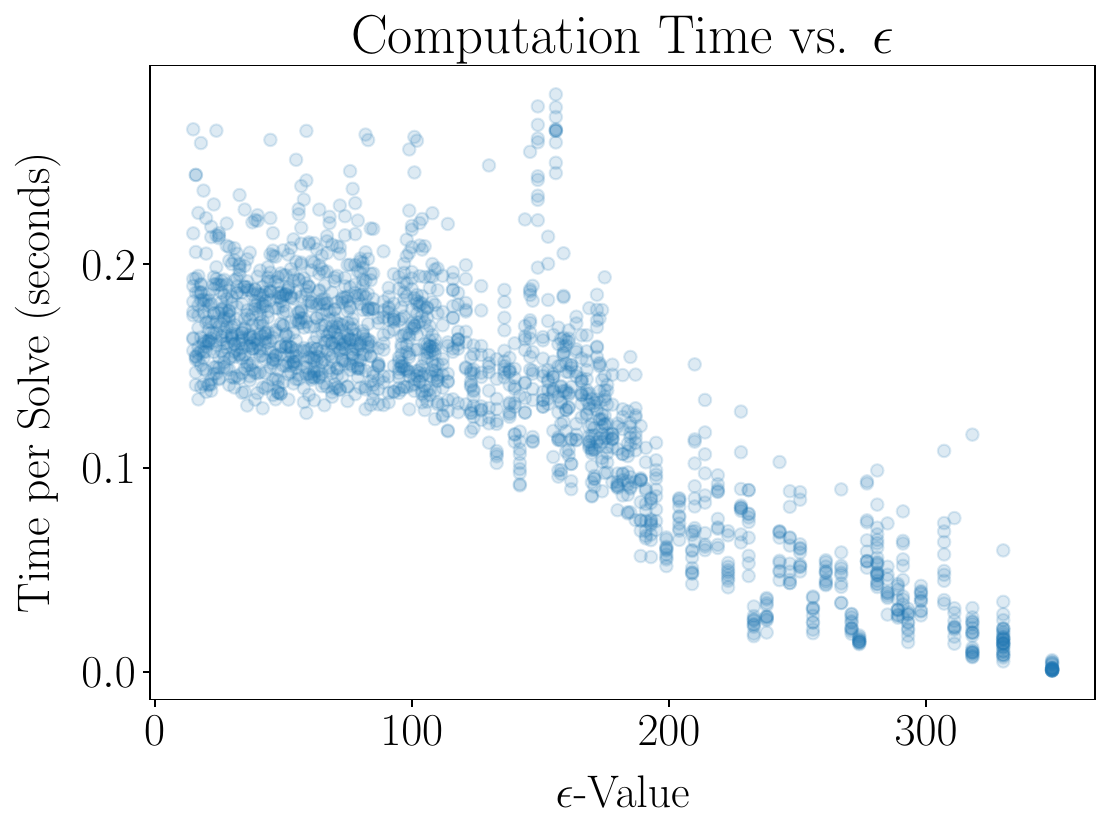}
  \end{minipage}\hfill
  \begin{minipage}[c]{0.19\textwidth}
    \caption{Computation time per solve for ILP plotted against $\epsilon$, for fixed sparsity $k = 20$.} \label{epsilon_sandwich.fig}
  \end{minipage}
\end{figure}

Problem \ref{ILP_problem} is especially easy to solve, for a fixed $k$, when $\epsilon$ is large, since constraints $\ref{ILPc1}$ and \ref{ILPc2} work together to sandwich down the combinatorial space of feasible solutions for $z$. We show this effect in Figure \ref{epsilon_sandwich.fig}, as $\epsilon$ increases the time it takes for an off-the-shelf optimization solver (Gurobi) to solve ILP Problem \ref{ILP_problem} drastically decreases. As such, Algorithm \ref{OA.alg} is especially efficient when $\epsilon$ is large.

\subsubsection{Efficient Gradient and Objective Evaluation:} In each iteration of Algorithm \ref{OA.alg}, we must evaluate objective $H_2(z)$ and compute gradient $\nabla H_2(z)$ for solution $z$. We show here that these computations can be conducted efficiently. Given solution vector $z$, let $\kappa$ represent the set of nonzero indices in $z$. The cardinality of set $\kappa$ is restricted to be at most $k$, the number of decision rules constructed in the rule set.  Furthermore, let $M_\kappa$ represent the sub-matrix of $M$ with the $\kappa$-indexed columns selected. Computing gradient  $\nabla H_2(z)$ and evaluating objective $H_2(z)$ is bottle-necked by the matrix inversion $(\frac{\mathbbm{I}}{\gamma} + M_\kappa^\intercal M_\kappa)^{-1}$, i.e., the cost of inverting a square matrix of at most size $k \times k$ \citep{bertsimas2020sparse}. Since we are interested in using MOSS to construct a sparse interpretable rule set, of no more than 10-20 decision rules, $k$ is typically small and this computation is cheap.

Algorithm \ref{OA.alg}, with the attributes discussed above, allows us to solve Problem \ref{mainoptimizationproblem} efficiently for a fixed value of $\epsilon$; we show timing results in \S\ref{computation_time_eval.section}. The algorithm is especially efficient when $\epsilon$ is large and when $k$ is small; $k$ is typically kept small in MOSS for interpretability reasons. However, we must solve Problem \ref{mainoptimizationproblem} across a range of $\epsilon$ values to compute the Pareto frontier between accuracy and stability. We discuss how to do so efficiently below.

\subsection{ Computing the Pareto Frontier}\label{efficient_pareto.section}

In this section, we introduce our method to efficiently compute the Pareto frontier between stability objective $H_1(z)$ and accuracy objective $H_2(z)$. To compute this Pareto frontier, we must repeatedly solve Problem \ref{mainoptimizationproblem} across a range of $\epsilon$ values, and we first discuss what values of $\epsilon$ to solve for below.

\subsubsection{$\epsilon$-Range:} We demonstrate here that it is sufficient to solve Problem \ref{mainoptimizationproblem} for a discrete sequence of $\epsilon$ values to fully compute the Pareto frontier between $H_1(z)$ and $H_2(z)$ with complete granularity. First, we note that for fixed sparsity $k$, the largest value of $\epsilon$ such that Problem \ref{mainoptimizationproblem} is feasible is equivalent to the $k$ largest elements of $\Pi$. We denote this value as $\epsilon_{\max}$ and recall that from \S\ref{knapsack.section} that $\epsilon_{\max}$ is equivalent to the optimal objective value of our stability selection optimization formulation (Problem \ref{knapsack.problem}).

To compute the Pareto frontier, we are interested in finding a sequence of $\epsilon$-values that correspond to a sequence of increasingly suboptimal objectives for $H_1(z)$. As discussed in \S\ref{knapsack.section}, this can be easily obtained by sorting $\Pi$ in descending order and taking the sum of a rolling window of $k$ elements down the vector to obtain objective values $H_1(z^*) > H_1(z^1) > \ldots > H_1(z^t)$, where $t = m - k + 1$. We set these objective values as the sequence $\mathbbm{E} = \epsilon_{max} > \epsilon^{1} > \ldots > \epsilon^t$. To compute the entire Pareto frontier between $H_1(z)$ and $H_2(z)$, in full granularity, we solve Problem \ref{mainoptimizationproblem} for every value of $\epsilon \in \mathbbm{E}$. It is important to note that we do not necessarily need to compute the Pareto frontier between $H_1(z)$ and $H_2(z)$ in full granularity. Rather, we can often compute sub-sequence or interesting segments of $\mathbbm{E}$ to find an appropriate model.

\subsubsection{Efficient Pareto Method:} Given a descending sequence of $\epsilon$-values $\mathbbm{E}$, we develop the following algorithm to efficiently solve Problem \ref{mainoptimizationproblem} for each value of $\epsilon \in \mathbbm{E}$. Our algorithm relies on the fact that we can exploit nested problem structure to reuse cutting planes, which greatly reduces the number of iterations needed for Algorithm \ref{OA.alg} to converge.

Say that we solve Problem \ref{mainoptimizationproblem} for two values of $\epsilon \in \mathbbm{E}$, where $\epsilon_1 > \epsilon_2$. Any solution to Problem \ref{mainoptimizationproblem} with $\epsilon = \epsilon_1$ is also a feasible solution to Problem \ref{mainoptimizationproblem} with $\epsilon = \epsilon_2$. As such, the cutting planes and approximation function generated when solving the $\epsilon_1$ problem are also valid for the $\epsilon_2$ problem. If we start with $\epsilon_{\max}$ and solve down the sequence $\mathbbm{E}$, we are solving a sequence of nested optimization problems, where every feasible solution found is also feasible for the subsequent problem. Consequently, each time we solve Problem \ref{mainoptimizationproblem} using Algorithm \ref{OA.alg}, we can re-use all of the cutting planes from earlier problems. This greatly reduces the number of iterations of Algorithm \ref{OA.alg} for each solve in the sequence. We present this efficient Pareto method (EPM) in the algorithm below.

\begin{algorithm}\label{EPM.alg}
    \caption{Efficient Pareto Method (EPM)}
    \footnotesize
    \textbf{Given:} $\mathbbm{E}$, $k$

    $Z \leftarrow \emptyset$ \tcp{Set of solutions.}
    
    $C \leftarrow \emptyset$ \tcp{Set of cutting planes.}

    $t = 0$, \ $z_t = \emptyset$, $c_t = \emptyset$ \tcp{Solution and cutting planes at iteration $t$.}
    
    \While{$\epsilon \in \mathbbm{E}$}{
        
        Warm start Algorithm \ref{OA.alg} with constraint set $C$ and current solution $z_t$.

        Apply Algorithm \ref{OA.alg} for $\epsilon$ and $k$ to obtain solution $z_{t+1}$ and cutting planes $c_{t+1}$.

        $Z \leftarrow Z \bigcup z_t$
        
        $C \leftarrow C \bigcup c_t$
        
        $t = t+1$}
    \textbf{return} $Z$
  \end{algorithm}
In practice, EPM greatly reduces the number of iterations required each time we apply Algorithm \ref{OA.alg} when we compute the Pareto frontier. We show this though an example below. 

\begin{figure}[h]
  \begin{minipage}[c]{0.28\textwidth}
    \includegraphics[width=\textwidth]{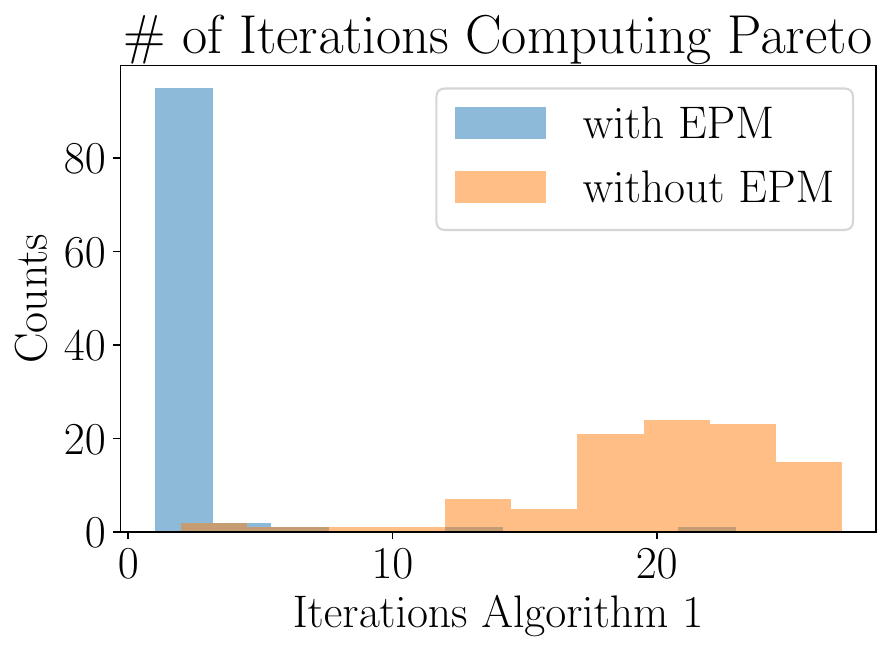}
  \end{minipage}\hfill
  \begin{minipage}[c]{0.18\textwidth}
    \caption{Our efficient EPM procedure reduces the number of OA iterations needed to compute the Pareto frontier.} \label{EPM_ablate_iters.fig}
  \end{minipage}
\end{figure}

We use MOSS with and without EPM to compute the first 100 values of $\mathbbm{E}$. In Figure \ref{EPM_ablate_iters.fig}, we show the number of iterations required each time we apply our cutting plane Algorithm \ref{OA.alg}. This corresponds to the number of cutting planes added, and the number of times we need to solve ILP Problem \ref{ILP_problem}. We see that EPM greatly reduces the number of iterations required, fact, often times only a few additional cutting planes are needed to solve Problem \ref{mainoptimizationproblem} for the next value of $\epsilon$. In \S\ref{computation_time_eval.section}, we show that this translates to large computational savings.

\subsection{Putting Together the Pieces} \label{putting_together_pieces.section}

By combining our new multi-objective optimization formulation in MOSS with our tailored cutting plane algorithm, we can efficiently explore the Pareto frontier between in-sample stability $H_1(z)$ and in-sample loss $H_2(z)$. Recall that we maximize in-sample stability as a proxy for empirical out-of-sample stability and minimize in-sample loss as a proxy for out-of-sample error. In this section, we apply MOSS to an example to explore how solutions along the Pareto frontier of $H_1(z)$ and $H_2(z)$ perform in terms empirical stability and out-of-sample $1 - R^2$.

On the \textsc{Galaxy} dataset from OpenML \cite{OpenML2013}, we first apply the SIRUS framework to generate a set of $m = 1000$ unique decision rules with selection proportions $\Pi$. We fix the sparsity of our desired rule set at $k = 15$ and apply MOSS to compute the entire Pareto frontier between in-sample stability and in-sample accuracy.

\begin{figure}[h]
    \centering
    \includegraphics[width = 0.46\textwidth]{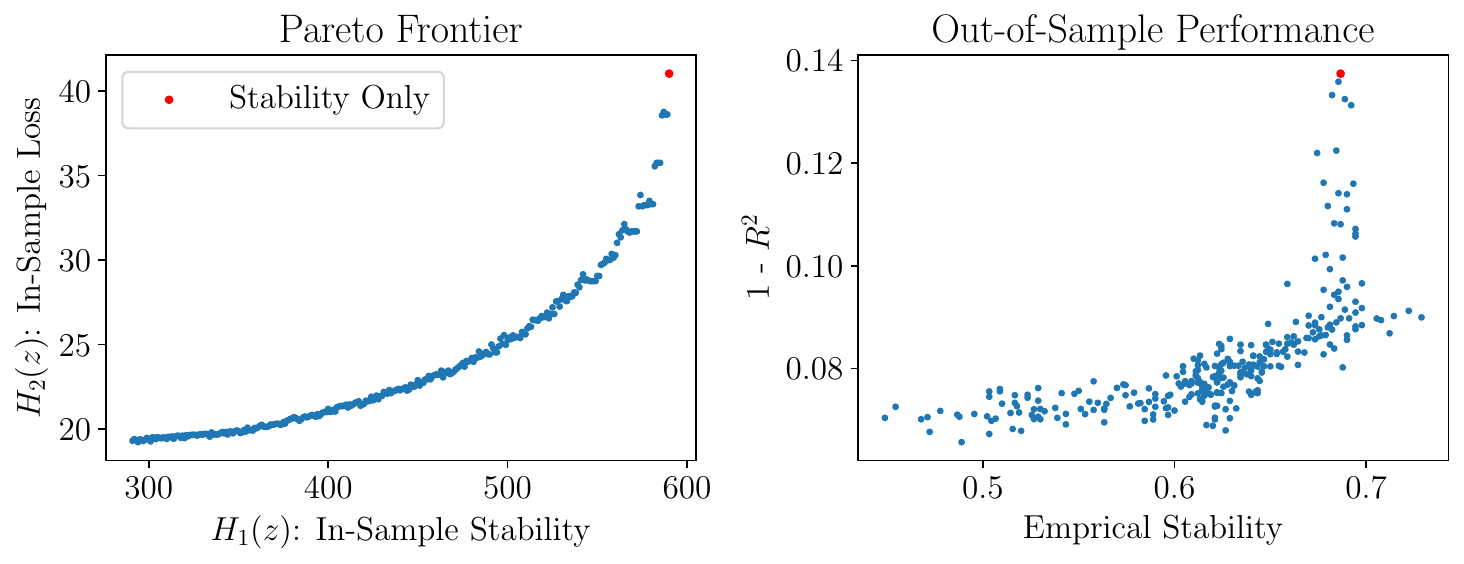}
    \caption{Left Panel: Pareto frontier between accuracy and stability recovered by MOSS. Note the steepness of the Pareto front near the stability-only solution. Right Panel: Out-of-sample accuracy and stability for MOSS solutions. MOSS can find stable solutions with much improved test accuracy. }
    \label{pareto_example.fig}
\end{figure}

The left panel in Figure \ref{pareto_example.fig} displays the Pareto frontier. The horizontal axis shows in-sample stability and the vertical axis shows in-sample loss, and the points show solutions of MOSS for different values of $\epsilon$. The red point at the top right of the Pareto indicates the solution when $\epsilon$ is set to its maximum value, meaning MOSS optimizes solely for in-sample stability. Note that the Pareto frontier drops sharply directly to the left of the red point. This is important; many solutions on the Pareto frontier slightly suboptimal in terms of in-sample stability have much lower in-sample loss.

Importantly, these solutions demonstrate better out-of-sample predictive performance while preserving good out-of-sample empirical stability. In the right panel in Figure \ref{pareto_example.fig}, we repeat the procedure detailed above across a 10-fold CV of the \textsc{Galaxy} dataset. On the horizontal axis, we report out-of-sample empirical stability, defined via the average DSC metric presented in \S\ref{model_stability_def.section}, and on the vertical axis we report out-of-sample performance, via 1-$R^2$, averaged across all folds. Again, each point shows MOSS solutions for different values of $\epsilon$. We observe that along the right hand side of the plot, many solutions have high empirical stability but some have much better out-of-sample performance; these correspond to the solutions that are slightly suboptimal in terms of in-sample stability discussed in the paragraph above.

We hypothesize that solutions in MOSS that are slightly suboptimal in terms of in-sample stability may fit the data much better, and produce models with much better out-of-sample performance. Moreover, these solutions may not necessarily exhibit worse out-of-sample stability, as we see above. Finally, as discussed in \S\ref{optalg.sec}, the cutting plane algorithm in MOSS is extremely efficient at exploring these solutions. We test this hypothesis in \S\ref{experiments.section} by evaluating how well MOSS performs compared to several state-of-the-art rule-based algorithms in terms of both accuracy and stability.

\subsection{Approximate Algorithm} 

As an aside, we also present an approximate algorithm to find high-quality solutions to Problem \ref{mainoptimizationproblem}. Our algorithm is self-contained and is much more computationally efficient than our outer-approximation approach when $\epsilon$ is small.

We first re-express Problem \ref{mainoptimizationproblem} in terms of just $w$, where $z_i = \mathbbm{I}(w_i \neq 0)$, and work with the unconstrained Lagrangian of this problem:
\begin{equation*}
    \min_{w}  \frac{1}{2}||y - Mw||_2^2 + \frac{1}{2\gamma}||w||_2^2 + \lambda_1 \sum_{i=1}^m \mathbbm{1}(w_i \neq 0) - \lambda_2 \sum_{i=1}^m \Pi_i \mathbbm{1}(w_i \neq 0).
\end{equation*}
We drop the constraints that involve $k$ and $\epsilon$ and use non-negative hyperparameters $\lambda_1$ and $\lambda_2$ to control the sparsity and in-sample stability of the extracted rule sets. This expression simplifies to:
\begin{equation}\label{heuristic.problem}
   \min_{w} \frac{1}{2} \biggl(||y - Mw||_2^2 + \frac{1}{\gamma}||w||_2^2 \biggr) + \biggl( \sum_{i=1}^m \mathbbm{1}(w_i \neq 0) (\lambda_1 - \Pi_i \lambda_2)\biggr),
\end{equation}
where the first term is smooth and the second term is separable over $w$. As such, we apply coordinate descent to Problem \ref{heuristic.problem} to find good solutions, since the algorithm will converge to a local minima. Each coordinate update has a closed-form solution, and we present our full coordinate descent heuristic in the appendix (\ref{CD.appx}).

\begin{figure}[h]
    \centering
    \includegraphics[width=0.47\textwidth]{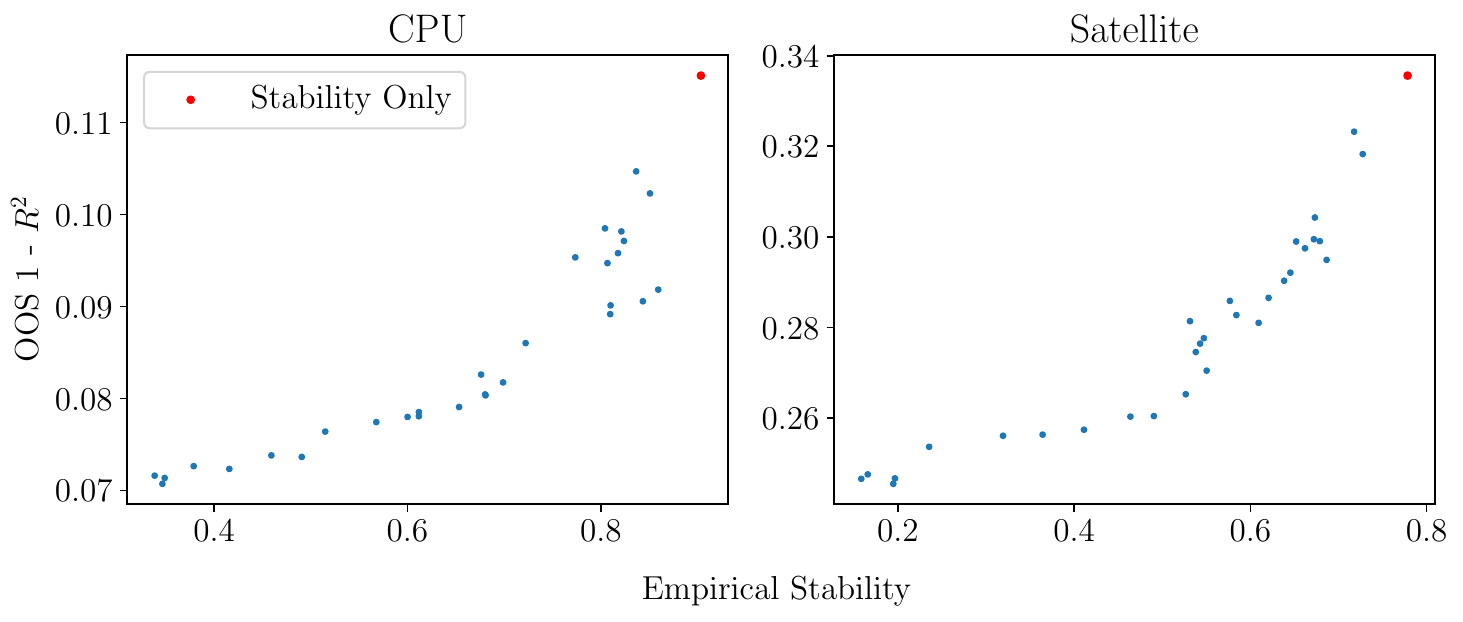}
    \caption{Our approximate algorithm recovers good approximations of the Pareto frontier between accuracy and stability.}
    \label{heuristic_paths.fig}
\end{figure}

In practice, we observe that our approximate algorithm can efficiently compute high-quality approximations of the Pareto frontier between accuracy and stability. In Figure \ref{heuristic_paths.fig}, we use our approximate algorithm to construct rule sets of sparsity $k = 15$ on the CPU and Satellite datasets from OpenML \cite{OpenML2013}. We sweep though parameter $\lambda_2$ to explore the trade-off between accuracy and stability, and we plot out-of-sample $1-R^2$ against empirical stability. We see from this figure that our approximate algorithm can find good solutions that balance out-of-sample accuracy with empirical stability.

\section{Experiments}
We evaluate the computation time and performance of MOSS.

\subsection{Computation Time Experiment} \label{computation_time_eval.section}

We first evaluate our cutting plane algorithm (Algorithm \ref{OA.alg}) against two state-of-the-art commercial solvers, Gurobi and MOSEK \cite{mosek,gurobi}. Using each method we solve various instances of Problem \ref{mainoptimizationproblem} across 10 fixed values of $\epsilon$, for rule sparsity $k = 15$. We conduct this timing experiment on a 2022 M2 Macbook Pro.

We show these computation time results in Table \ref{timing.table}. The leftmost column shows the problem size in terms of the number of data points $n$ and the number of decision variables $m$. We observe that on the smallest problems our cutting plane algorithm achieves orders of magnitude speedups compared to commercial solvers. For larger problems, we  observe that our cutting plane algorithm scales well, and can handle problem sizes beyond the capabilities of off-the-shelf optimization solvers.

\begin{table}[h]
\scalebox{0.8}{
\begin{tabular}{|c|c|c|c|}
\hline
\textbf{Data Points / Variables} & \textbf{Cutting Plane} & \textbf{Gurobi}                               & \textbf{MOSEK}                                \\ \hline
\textbf{150, 250}                & 0.044s (0.005)               & \multicolumn{1}{l|}{28m 15s}                  & \multicolumn{1}{l|}{40m 12s}                  \\ \hline
\textbf{150, 500}                & 0.0965s (0.003)              & 31m 10s                                       & 45m 6s                                        \\ \hline
\textbf{150, 1000}               & 0.464s (0.68)                & \multicolumn{1}{l|}{1h 30m}                   & \multicolumn{1}{l|}{1h 42m}                   \\ \hline
\textbf{1500, 1200}              & 0.66s (0.09)                 & \cellcolor[HTML]{FFCCC9}                      & \cellcolor[HTML]{FFCCC9}                      \\ \hline
\textbf{5000, 500}               & 3.07s (0.1)                  & \cellcolor[HTML]{FFCCC9}                      & \cellcolor[HTML]{FFCCC9}                      \\ \hline
\textbf{7500, 1500}              & 2.27s (0.9)                  & \multicolumn{1}{l|}{\cellcolor[HTML]{FFCCC9}} & \multicolumn{1}{l|}{\cellcolor[HTML]{FFCCC9}} \\ \hline
\textbf{7500, 2500}              & 2.3s (0.5)                   & \multicolumn{1}{l|}{\cellcolor[HTML]{FFCCC9}} & \multicolumn{1}{l|}{\cellcolor[HTML]{FFCCC9}} \\ \hline
\textbf{15000, 700}              & 14.1s (0.3)                  & \cellcolor[HTML]{FFCCC9}                      & \cellcolor[HTML]{FFCCC9}                      \\ \hline
\end{tabular}}
\caption{Computation time results, the red cells indicate that the method fails to return the optimal solution after 4 hours.}
\label{timing.table}
\end{table}

\subsubsection{Ablation Study} We also investigate the impact of our efficient Pareto method on computation time. Our efficient Pareto method exploits nested problem structure to reuse cutting planes are warm-starts when computing the Pareto frontier, so we compare this method against naively sweeping through values of $\epsilon$.

\begin{table}[h]
\scalebox{0.85}{
\begin{tabular}{|c|c|c|}
\hline
\textbf{Data Points / Variables} & \textbf{w/ EPM} & \textbf{w/o EPM} \\ \hline
\textbf{1000, 500}               & 13s (2.1)       & 26s (1.4)        \\ \hline
\textbf{3000, 600}               & 65s (5.2)       & 5m 10s (24.1)    \\ \hline
\textbf{5000, 500}               & 58s (4.3)       & 9m 12s ( 32.1)   \\ \hline
\textbf{8000, 400}               & 9m 5s (10.2)    & 55m 12s (42.3)   \\ \hline
\end{tabular}} 
\caption{Timing results for ablation study.}\label{ablation_study.table}
\end{table}
\vspace*{-\baselineskip}

On several instance of Problem \ref{mainoptimizationproblem} of varying sizes, we use MOSS to compute the Pareto frontier across 100 values of $\epsilon$, with and without our efficient Pareto method (EPM). We show the results in Table \ref{ablation_study.table}. We see from this table that EPM drastically reduces the computation time required for MOSS to compute Pareto frontiers and can yield up to an order of magnitude speedup for larger problem sizes.

\subsection{Performance Experiment} \label{experiments.section}
We evaluate MOSS against competing rule-based algorithms in terms of stability and predictive performance.

\subsubsection*{\textbf{Experiment Procedure}} We repeat this procedure on 30 regression datasets of various sizes sourced from OpenML \cite{OpenML2013}; the full lists of datasets with metada can be found in the appendix (\ref{dataused.appx}). On each dataset, we conduct a 10-fold cross validation and on each fold we fit a random forest to generate $m \sim 10^3$ unique candidate rules with selection proportions $\Pi$. We apply MOSS to construct sets of $k = 15$ decision rules under three settings. For the high $\epsilon$ setting (\textbf{MOSS-$\epsilon$-H}) we set $\epsilon$ to be the 3rd element in sequence $\mathbbm{E}$, close to $\epsilon_{\max}$, and for the medium setting (\textbf{MOSS-$\epsilon$-M}) we set $\epsilon$ to be the 40th element of $\mathbbm{E}$, close to the midpoint of the sequence. For these two settings, we compute solutions using our OA algorithm with EPM. For the low setting (\textbf{MOSS-$\epsilon$-L}) we compute solutions using our approximate algorithm. Across all folds, we evaluate the out-of-sample empirical stability, defined via average pairwise DSC, and predictive performance, defined via average out-of-sample $R^2$, of the constructed rule sets.

We compare MOSS against the following state-of-the-art competing algorithms: \textbf{SIRUS} (2021) \cite{benard2021sirus} which constructs decision rules with respect to stability only, \textbf{FIRE} (2023) \cite{liu2023fire} which selects rule sets using the non-convex MCP penalty,  \textbf{GLRM} (2019)\cite{wei2019generalized}  which constructs compact rule sets from scratch using column generation, and \textbf{RuleFit} (2008) \cite{friedman2008predictive} which selects rules using the LASSO. For a fair comparison, we use each algorithm to construct rule sets of $k=15$ decision rules. Again, we assess the out-of-sample empirical stability and predictive performance of all algorithms across all folds. Further details about our experiment can be found in the appendix (\ref{experimentaldetails.appx}).

\begin{figure}[h]
    \centering
    \includegraphics[width=0.38\textwidth]{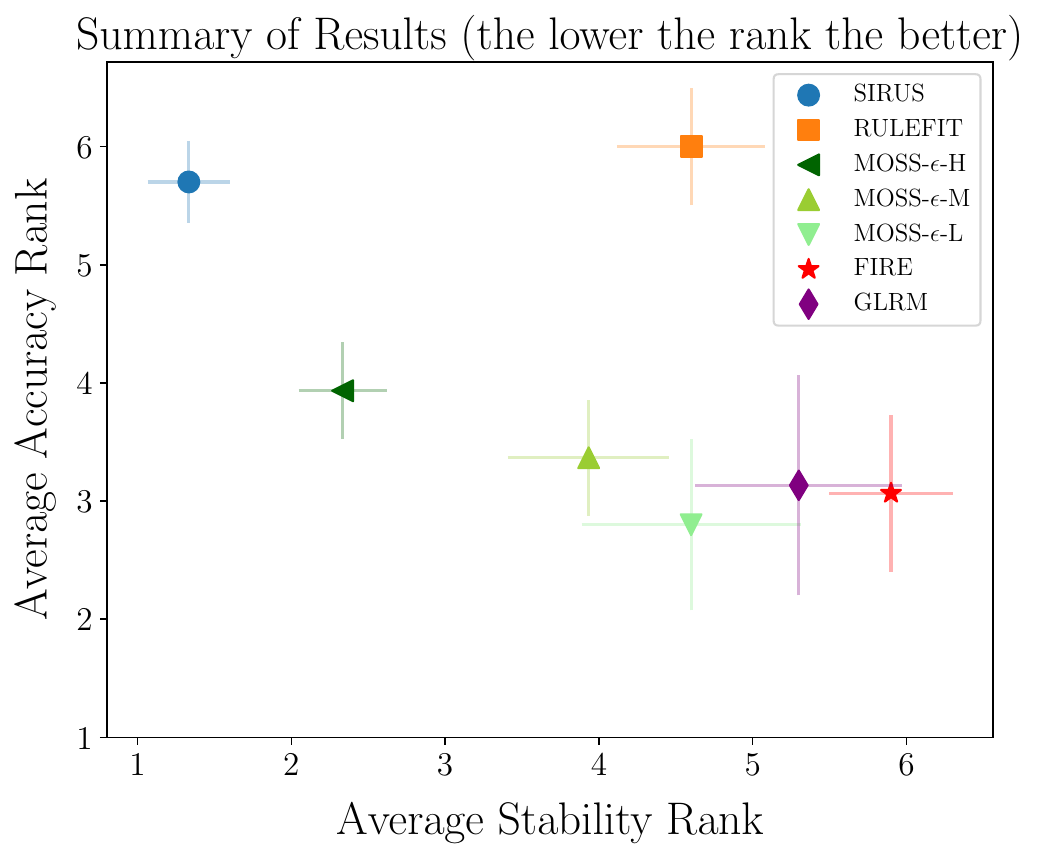}
    \caption{Method rank in terms of accuracy and stability.}
    \label{method_ranks.fig}
\end{figure}

\subsubsection*{\textbf{Results}} We discuss our experimental results; tables of the detailed results for each dataset can be found in the appendix (\ref{detailed_results.appx}). We first consider the ranking of each method, in terms of accuracy and predictive performance, on each dataset in the experiment (with 1 being the best and 7 being the worst). In Figure \ref{method_ranks.fig}, we plot the average rank of each method where the horizontal axis shows accuracy and the vertical axis shows stability. The points show which method is used and the lines show standard error. We observe from this figure that SIRUS performs the best in terms of model stability but nearly the worst in terms of model accuracy. Among the competing methods, \textsc{FIRE} and GLRM perform the best in terms of accuracy but the worst in terms of stability. RuleFit performs the worst in terms of accuracy but is more stable than \textsc{FIRE}, which may be due to the added shrinkage of the LASSO \cite{liu2023fire}.

Our MOSS methods are unique in that they perform well with respect to both objectives. MOSS-$\epsilon$-H is the second most stable method, and is much more accurate compared to SIRUS and RuleFit. Interestingly, MOSS-$\epsilon$-L performs the best in terms of model accuracy and is still more stable than \textsc{FIRE} and GLRM. This may be due to the fact that the combined sparsity-stability penalty in our approximate algorithm leads to better out-of-sample performance compared to penalizing sparsity only. As expected, MOSS-$\epsilon$-M balances stability and accuracy compared to the  competing methods. 

To condense accuracy and stability into a single metric we take the average combined rank of each method across all dataset. We show these results in the table below and it is apparent that MOSS-$\epsilon$-H performs the best in terms of this metric. As such, we recommend setting MOSS-$\epsilon$-H as the default when using our framework.

\begin{table}[h]
\scalebox{0.775}{
\begin{tabular}{|c|c|c|c|c|c|c|}
\hline
\cellcolor[HTML]{EFEFEF}\textbf{MOSS-$\epsilon$-H} & \textbf{MOSS-$\epsilon$-M} & \cellcolor[HTML]{EFEFEF}\textbf{MOSS-$\epsilon$-L} & \cellcolor[HTML]{EFEFEF}\textbf{SIRUS} & \cellcolor[HTML]{EFEFEF}\textbf{FIRE} & \cellcolor[HTML]{EFEFEF}\textbf{GLRM} & \cellcolor[HTML]{EFEFEF}\textbf{RULEFIT} \\ \hline
3.13                                               & 3.65                       & 3.74                                               & 3.51                                   & 4.18                                  & 4.51                                  & 5.3                                      \\ \hline
\end{tabular}}
\end{table}

With this in mind, we compare MOSS-$\epsilon$-H against SIRUS and \textsc{FIRE}, the competing algorithms that perform the best in terms of stability and accuracy individually. On 24 out of the 31 datasets in our experiment, the empirical stability of the SIRUS rule set is within one standard error of MOSS-$\epsilon$-H, making rule sets functionally identical in terms of stability. On these datasets, MOSS-$\epsilon$-H is much more accurate and exhibits an average \textbf{10\%} increase in out-of-sample $R^2$. Compared to $\textsc{FIRE}$, MOSS-$\epsilon$-H exchanges a \textbf{2\%} decrease in out-of-sample $R^2$ for a \textbf{190\%} increase in empirical stability. We show the distribution of these percent differences in the appendix (\ref{comparisons.appx}). As our experiments show, MOSS can construct rule sets that jointly out-perform our competing methods in terms of both accuracy and stability.

\section{Conclusion}

We conclude by showcasing the utility of MOSS on two real-world case studies. We also discuss connections between MOSS and related works and  analyze the sensitivity of our framework to parameters $\gamma$ and $k$. 

\subsubsection*{\textbf{Case Studies}} In Section \ref{case_study.appx} of the appendix, we present two case studies where we apply MOSS to real-world problems: scientific discovery in marine biology and census planning. In both of these examples, we show that MOSS can construct stable rule sets that perform well and reveal uncover insights about the data.

\subsubsection*{\textbf{Related Work}} We explore connections between MOSS and Rashomon set algorithms in interpretable machine learning \citep{xin2022exploring,mata2022computing,ciaperoni2024efficient}. Rashomon set algorithms identify a collection of models—such as decision rules—that achieve near optimal predictive performance in terms of training error. Despite their similar accuracy, these models can vary significantly in other aspects, such as fairness or the features they use. By analyzing this so-called Rashomon set of near optimal models, practitioners can account for these additional considerations when selecting an appropriate.

Consider the task of extracting rule sets from tree ensembles. Under this setting, we see an example of the Rashomon effect in Figure \ref{knapsack_explore.fig}. In this plot, the extracted rule sets indicated by the blue points labeled 2 through 10 have similar in-sample training errors, however, the in-sample stability of these models vary significantly.

With this in mind, we can consider an alternative formulation of MOSS that explores the Rashomon set of size-$k$ rule sets that can be extracted from a tree ensemble, and selects the rule set with the highest in-sample stability. This formulation is given by:
\begin{equation}\label{rashomonproblem}
\begin{alignedat}{2}
&\min_{w, z} &&\quad - H_1(z) , \\
&\text{s.t.} &&\frac{1}{2}\Vert y -  \sum_{i=1}^m f_i(X) w_i\Vert_2^2 \leq L^* + \psi, \quad \sum_{i=1}^m z_i \leq k, \\
&&& \quad w_i(1-z_i) = 0 \quad  \forall \ i \in [m],  \quad z \in \{0,1\}^m,
\end{alignedat}
\end{equation}
where $L^*$ is the lowest possible training error achieved by extracting a size-$k$ rule set from the tree ensemble and $\psi$ is the error tolerance \citep{xin2022exploring}. Problem \eqref{rashomonproblem} is similar to MOSS if we were to move  objective $H_2(z)$ into the constraint set for our $\epsilon$-constraint method. However, as we discuss in \S\ref{BIP.section}, moving the in-sample loss objective into the constraint set results in a much more difficult problem to solve due to its non-linearity. As such, MOSS explores the collection of size-$k$ rule sets with near optimal in-sample stabilities and selects the model with the lowest in-sample loss. This is similar to a Rashomon set method, but it leads to a form more amenable to optimization.  

It is important to note that while Rashomon set methods are interesting, enumerating or exploring the Rashomon set is extremely computationally challenging. As such, current Rashomon set methods for decision rules are restricted to binary classification tasks \citep{xin2022exploring,ciaperoni2024efficient} which prevents direct comparisons with our method. MOSS explores the Pareto frontier of just two objectives, accuracy and stability, extremely efficiently so that practitioners can rapidly select an appropriate model.

\subsubsection*{\textbf{Sensitivity Analyses}} In  Section \ref{sens.section} of the appendix,  we evaluate the sensitivity of MOSS to two parameters: $\gamma$, which controls the regularization penalty in the accuracy objective, and $k$, which controls the size of the constructed rule sets. In \ref{gammasense.section}, we demonstrate that MOSS is relatively \emph{insensitive} to $\gamma$ in terms of both accuracy and stability.  In \ref{ksens.section}, we show that increasing $k$ improves the accuracy and stability of rule sets, but reduces interpretability. Additionally, we show results for our experiment in \S\ref{experiments.section} for different values of $k$. Based on these findings, we recommend selecting $\gamma \in [10^{-3}, 10^{-2}]$ and $k \in \{10,15,20\}$ as good starting points for MOSS.

\subsubsection*{\textbf{Final Remarks}}  MOSS is a useful framework that quickly constructs rule sets to explore the Pareto frontier between stability and accuracy. By examining this trade-off, practitioners can select models suitable for drawing trustworthy insights from the data.

\subsubsection*{\textbf{Acknowledgments.}} The authors acknowledge support from the Office of Naval Research (ONR) under grants N000142212665 and N000142112841, and from the MIT Health Systems Initiative.


\bibliographystyle{plainnat}
\balance
\bibliography{conference/confv1.bib}

\newpage
\section*{Appendix}
\setcounter{equation}{0}
\renewcommand{\thesubsection}{\Alph{subsection}}

Our appendix can be found below. Code to reproduce our experiments can be found in this repository: \href{https://github.com/brianliu12437/MOSS}{github.com/brianliu12437/MOSS}

\subsection{Proofs}

We show the proofs for our propositions below.

\subsubsection{Proof of Proposition 1}
\label{prop1.appx}

First, given decision vector $z$, let $s$ be the set of indices such that $z \neq 0$. Let matrix $M_s$ correspond to the sub-matrix of $M$ with the $s$ indexed rows. We have that:
\begin{equation}
    \sum_{i = 2}^m z_i M_i M_i^{\intercal} = M_s (M_s)^\intercal.
\end{equation}

We start with function:
\begin{equation}
H_2(z) = \min_{w_1, \ldots, w_n}  \frac{1}{2}\Vert y - Mw\Vert_2^2 + \frac{1}{2\gamma} \Vert w \Vert_2^2    
\end{equation}
We have that minimizer $w^* = \biggl(\frac{\mathbbm{I}}{\gamma} - M^\intercal M \biggr)^{-1} M ^\intercal y$. Plugging in $w^*$ into $\frac{1}{2}\Vert y - Mw\Vert_2^2 + \frac{1}{2\gamma} \Vert w \Vert_2^2$  evaluating the objective, and combining with the expression above, yields the desired expression for $H_2(z)$.

\subsubsection{Proof of Proposition 2}
\label{prop2.appx}

We start with function:
\begin{equation*}
H_2(z) = \frac{1}{2}y^{\intercal}\biggl(\mathbb{I}_n +\gamma \sum_{i=1}^m z_i M_i M_i^{\intercal}\biggr)^{-1} y.   
\end{equation*}

Let matrix $A = \biggl(\mathbb{I}_n +\gamma \sum_{i=1}^m z_i M_i M_i^{\intercal}\biggr)$. We can re-express $H_2(z) = f(A) = y^\intercal A^{-1}y$. Matrix $A$ is symmetric positive definite
so function $f$ is convex over $A$ (follows from example 3.4 \cite{boyd2004convex}).

Define function $g(z) = \biggl(\mathbb{I}_n +\gamma \sum_{i=1}^m z_i M_i M_i^{\intercal}\biggr)$; function $g$ is affine over $z$. Convex combinations of affine functions are convex so $H_2(z) = f(g(z))$ is convex completing the proof.

\subsubsection{Convergence of Cutting Plane Algorithm} \label{alg_convergence.section}

The convergence of our cutting plane algorithm (Algorithm \ref{OA.alg}) is given by Corollary \ref{converge.corollary}.

\begin{corollary} \label{converge.corollary} Algorithm \ref{OA.alg} terminates after a finite number of cutting planes and returns $z^*$, the optimal solution to Problem \ref{mainoptimizationproblem}.
\end{corollary}

This follows directly from Fletcher and Leyffer (1994) \cite{fletcher1994solving}.

\subsection{Coordinate Descent Heuristic}\label{CD.appx}
In this section, we present our coordinate descent-based heuristic to find good solutions to optimization problems in MOSS.

Our CD heuristic finds good solutions to the problem:
\begin{equation*}\label{heuristic1.problem}
   \min_{w} \frac{1}{2} \biggl(||y - Mw||_2^2 + \frac{1}{\gamma}||w||_2^2 \biggr) + \biggl( \sum_{i=1}^m \mathbbm{1}(w_i \neq 0) (\lambda_1 - \Pi_i \lambda_2)\biggr).
\end{equation*}

We cycle through indices $k \in 1 \ldots m$ and update each decision variable $w_k$ one-at-a-time.

\textbf{Cyclic Updates:} Given a fixed index $k$, let $\delta$ represent the set of remaining indices $\{1 \ldots m\} \setminus k$. We define the residual vector as $r_k = y - \sum_{i \in \delta} M_i w_i$. Each update in our coordinate descent algorithm aims to solve the problem:
\begin{equation}\label{cd_update.problem}
    \min_{w_k} \frac{1}{2}||r - M_k w_k||_2^2 + \frac{1}{2\gamma}(w_k)^2 + \mathbbm{1}(w_k \neq 0) (\lambda_1 - \Pi_k \lambda_2).
\end{equation}
We can solve this problem efficiently by considering two scenarios. If $w_k = 0$, we have that the objective value of Problem \ref{cd_update.problem} is equal to $\frac{1}{2}||r||_2^2$. If $w_k \neq 0$, we have that:

\begin{align*}
    w_k^* = \frac{M_k^\intercal r}{M_k^\intercal M_k + \frac{1}{\gamma}}, 
\end{align*} and that the objective is equal to \begin{align*}
    \frac{1}{2}||r - M_k w_k^*||_2^2 + \frac{1}{2\gamma}(w_k^*)^2 + \lambda_1 - \Pi_k \lambda_2.
\end{align*}
  We take the scenario with the lower objective value as the solution for Problem \ref{cd_update.problem}.

  For our CD heuristic, we sweep through coordinates $k \in \{1\ldots m\}$ while solving Problem \ref{cd_update.problem}, and repeat sweeps until our algorithm converges. The algorithm will converge to a local minima since the first term in the objective is smooth and the second term is separable over $w$. In practice, we observe that it is much more efficient to apply our CD heuristic to find good solutions when $\lambda_2$ is small as opposed to the optimization problem to optimality for small values of $\epsilon$ using our cutting plane algorithm.

\subsection{Real World Case Studies} \label{case_study.appx}
Here, we present two case studies to showcase the usefulness of MOSS on real-world examples.

\subsubsection{Case Study: Stability in Scientific Discovery}
We apply MOSS to a real-world problem: drawing scientific conclusions from observational studies in ecology. Stability (replicability between repeated analyses) is critical to reliable scientific discovery. In Gagn\'e et al. 2018, researchers from the Monterey Bay Aquarium, US Fish and Wildlife Service, and multiple universities sample seabird feathers from multiple locations in the North Atlantic, North Pacific, and South Pacific oceans \cite{gagne2018seabird}. They analyze amino acid compounds in each sample to determine the estimated trophic position of the seabird, i.e. the relative position of the animal in its food chain. Organism with higher trophic levels consume organisms with lower ones. The researchers assemble a dataset of 18000 samples and 11 covariates that capture information about the seabird and its surrounding ecosystem. The goal of their study is to understand how these features contribute to trophic position.

To accomplish this, the researchers fit a random forest model to predict trophic position; the model performs well, with an out-of-sample $R^2$ of \textbf{0.58}. The authors examine the feature importance rankings of the model, determined via in-built variable importance scores \cite{breiman2001random}, and conduct a sensitivity analysis to assess the stability of these rankings by perturbing and re-analyzing the data. Bird species and human fishing pressure (catch per unit area) are consistently identified as important drivers of trophic position. This serves as the starting point for our application of MOSS.

\begin{figure}[h]
    \centering
    \includegraphics[width=.9\linewidth]{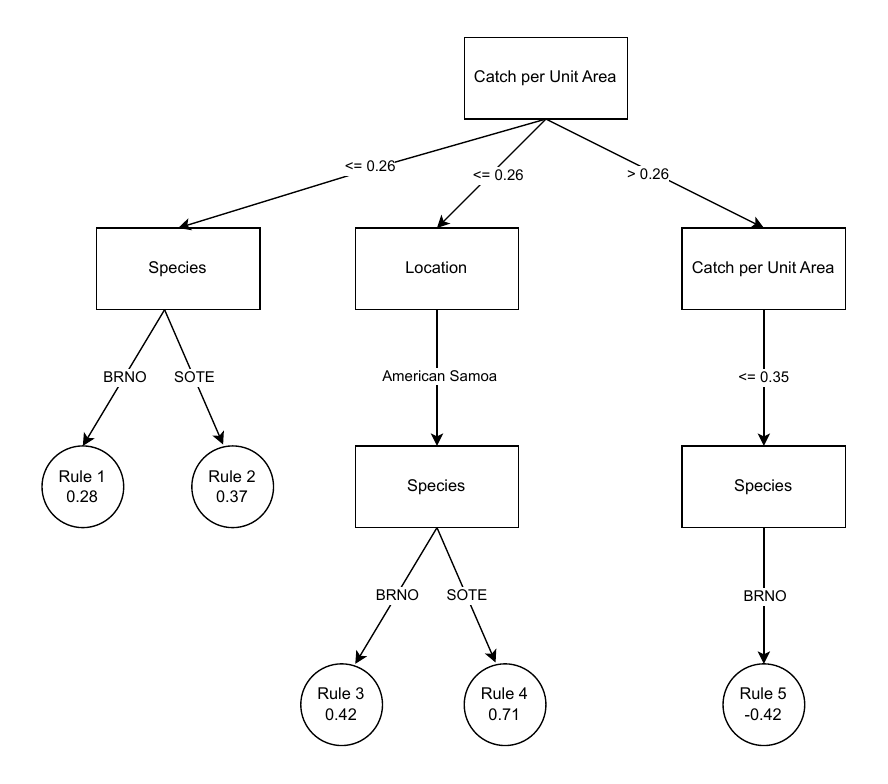}
    \caption{Stable structure of 5 decision rules consistently constructed by MOSS across analyses for the seabird case study.}
    \label{consistent_structure.fig}
\end{figure}

Our goal is to construct a sparse subset of decision rules that fit the data well and that are stable across re-analyses. These decision rules will provide us with a more granular understanding of how the covariates impact trophic position, compared to feature importance rankings. When we apply SIRUS to construct 10 decision rules, we obtain an out-of-sample $R^2$ of \textbf{0.30}, nearly half that of the original random forest. This rule set is not accurate enough to generate reliable conclusions. When we apply FIRE, we obtain an out-of-sample $R^2$ of \textbf{0.59}, however, \textbf{0} rules are consistently constructed across repeated analyses. These results are not stable, or replicable, enough to be considered trustworthy. When we apply MOSS, we obtain an out-of-sample $R^2$ of \textbf{0.57}, comparable to that of the random forest and we can identify a structure of \textbf{5} decision rules are consistently constructed across re-analyses. We show this structure in Figure \ref{consistent_structure.fig}.

By examining this structure, we can determine that a moderate increase in fishing pressure (catch per unit area) decreases the trophic position of the brown noody (BRNO), that seabirds in America Samoa have higher trophic positions, and the the sooty tern (SOTE) generally has higher trophic positions than the brown noody. These findings are more granular than the feature importance rankings presented in \cite{gagne2018seabird}, and are stable across repeated analyses.

\subsubsection{Case Study: Census Planning}

We apply MOSS to the real-world policy problem of identifying how demographic features at the census tract level contribute to low response rates for the American Community Survey (ACS). The ACS is administered annually by the U.S. Census Bureau and collects detailed demographic information on a sample of the U.S. population. ACS data is frequently used to inform policy decisions, including environmental justice and climate action initiatives \citep{lee2019census}. As such, ensuring that the data is collected from a representative sample of the population is important for policymaking. 

We use data from the U.S. Census Bureau Census Planning Database \citep{us_census_bureau_planning_database} to predict response rates for the 2017 ACS in California. Our dataset contains demographic information at the census tract level and consists of 8000 observations and 3000 features. Motivated by the observation that tree ensembles have historically performed well in predicting census return rates \citep{erdman2017low}, we first fit a random forest model. This black-box ensemble performs well, achieving an out-of-sample $R^2$ score of \textbf{0.79}, but lacks transparency.

To improve transparency, we construct a compact collection of decision rules. FIRE generates rule sets of 15 rules that achieve an out-of-sample $R^2$ score of \textbf{0.75}; however, these rules are highly unstable across repeated analyses, making them unreliable. To enhance stability, we apply SIRUS to construct rule sets; but this leads to a decline in predictive performance, reducing the out-of-sample $R^2$ score to \textbf{0.52}. 

\begin{figure}[h]
    \centering
    \includegraphics[width=.9\linewidth]{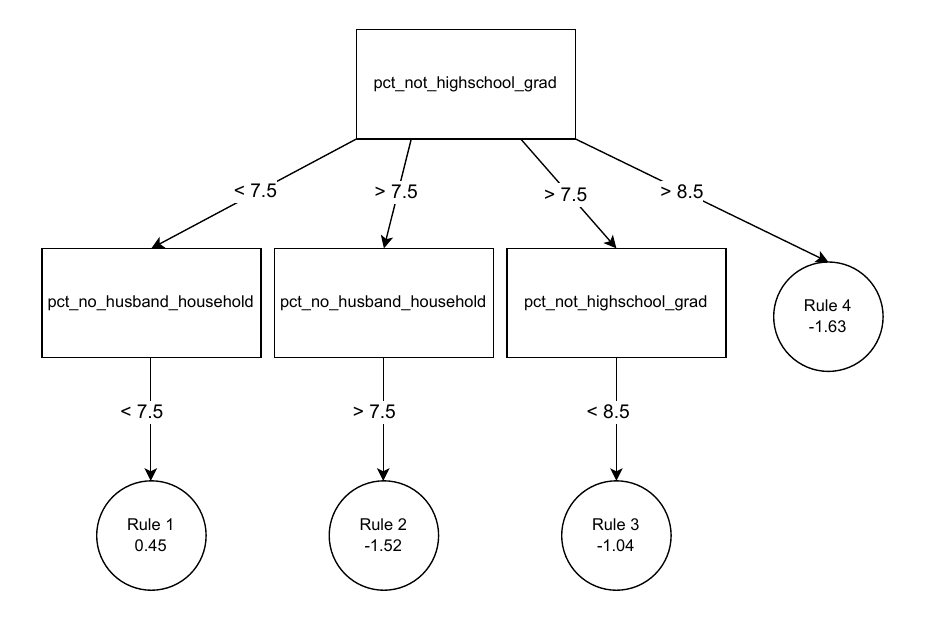}
    \caption{Stable structure of 4 decision rules consistently constructed by MOSS across analyses for the census planning case study.}
    \label{census_stab.fig}
\end{figure}

When we apply MOSS to construct a compact rule set, we achieve an out-of-sample $R^2$ score of \textbf{0.73}, significantly higher than that of SIRUS. Additionally, we identify a set of four decision rules that are consistently reproduced across re-analyses, which we show in Figure \ref{census_stab.fig}. From this structure, we observe that high school educational attainment and the percentage of households with no husband present are consistently identified as key drivers of low response rates. These findings may help inform strategies to improve response rates for the ACS.

\clearpage
\newpage

\newpage

\subsection{Experiment Details}
In this section, we discuss additional details regarding our experiments.

\subsubsection{Datasets Used} \label{dataused.appx}

We show a full table of the datasets used in our experiment in \S\ref{experiments.section}, along with the size of each dataset. All of these datasets were sourced from the OpenML repository \cite{OpenML2013}.

\begin{table}[h] \scalebox{0.98}{
\begin{tabular}{|c|c|c|c|}
\hline
\textbf{OpenML ID} & \textbf{Name}       & \textbf{Observations} & \textbf{Features} \\ \hline
296                & Ailerons            & 13750                 & 40                \\ \hline
1027               & ESL                 & 488                   & 4                 \\ \hline
42570              & mercedes            & 4209                  & 376               \\ \hline
41021              & Moneyball           & 1232                  & 14                \\ \hline
183                & abalone             & 4177                  & 8                 \\ \hline
196                & autoMpg             & 398                   & 7                 \\ \hline
195                & auto\_price         & 159                   & 15                \\ \hline
558                & bank32nh            & 8192                  & 32                \\ \hline
560                & bodyfat             & 252                   & 14                \\ \hline
227                & cpu\_small          & 8192                  & 12                \\ \hline
216                & elevators           & 16599                 & 18                \\ \hline
574                & house\_16H          & 22784                 & 16                \\ \hline
537                & houses              & 20640                 & 8                 \\ \hline
189                & kin8nm              & 8192                  & 8                 \\ \hline
405                & mtp                 & 4450                  & 202               \\ \hline
344                & mv                  & 40768                 & 10                \\ \hline
547                & no2                 & 500                   & 7                 \\ \hline
201                & pol                 & 15000                 & 48                \\ \hline
529                & pollen              & 3848                  & 4                 \\ \hline
308                & puma32H             & 8192                  & 32                \\ \hline
294                & satellite\_image    & 6435                  & 36                \\ \hline
541                & socmob              & 1156                  & 5                 \\ \hline
507                & space\_ga           & 3107                  & 6                 \\ \hline
223                & stock               & 950                   & 9                 \\ \hline
505                & tecator             & 240                   & 124               \\ \hline
315                & us\_crime           & 1994                  & 127               \\ \hline
519                & vinnie              & 380                   & 2                 \\ \hline
690                & visualizing\_galaxy & 323                   & 4                 \\ \hline
503                & wind                & 6574                  & 14                \\ \hline
287                & wine\_quality       & 6497                  & 11                \\ \hline
\end{tabular}}
\end{table}

\subsubsection{Experimental Details}\label{experimentaldetails.appx}

We conduct all of our performance experiments on a 2022 MacBook Pro. We implement MOSS in Python, and we will open-source our implementation after the review period. We will also open-source the code to reproduce all of our experiments. For our competing algorithms, we use the following implementations.
\begin{itemize}
    \item FIRE: We use the open-source Python implementation of FIRE found in the GitHub repository linked in \cite{liu2023fire}.
    \item GLRM: This method is implemented in the \textsc{AIX360} package maintained by IBM \cite{aix360-sept-2019}
    \item RuleFit: We implement the RuleFit algorithm in \textsc{scikit-learn} \cite{sklearn_api}.

    \item We implement the SIRUS algorithm found in \cite{benard2021sirus} in Python.
\end{itemize}

\clearpage
\newpage

\subsection{Experiment Results: Detailed Results} \label{detailed_results.appx}
In the pages below, we show tables (\ref{performance_means.table},\ref{performance_stds.table},\ref{stability_means.table},\ref{stability_std.table}) of our detailed experiment results by dataset. We report out-of-sample $R^2$ to measure predictive accuracy and average pairwise DSC to measure empirical stability. We report both the means and standard errors of both metrics over the 10-fold CV.

\begin{table*}[h]
\begin{tabular}{|c|c|c|c|c|c|c|c|}
\hline
\textbf{Dataset}         & \textbf{MOSS-$\epsilon$-H} & \textbf{MOSS-$\epsilon$-M} & \textbf{MOSS-$\epsilon$-L} & \textbf{FIRE} & \textbf{GLRM} & \textbf{SIRUS} & \textbf{RULEFIT} \\ \hline
Ailerons\_296            & 0.639792                   & 0.656669                   & 0.670918                   & 0.675325      & 0.324151      & 0.569526       & 0.546867         \\ \hline
ESL\_1027                & 0.801289                   & 0.828360                   & 0.801771                   & 0.814266      & 0.845445      & 0.756764       & 0.672672         \\ \hline
mercedes\_42570          & 0.572542                   & 0.556255                   & 0.576232                   & 0.574502      & 0.412054      & 0.570803       & 0.566909         \\ \hline
Moneyball\_41021         & 0.869957                   & 0.865513                   & 0.871525                   & 0.871419      & 0.928186      & 0.856857       & 0.683138         \\ \hline
abalone\_183             & 0.376287                   & 0.393225                   & 0.447597                   & 0.453783      & 0.508449      & 0.360119       & 0.348036         \\ \hline
autoMpg\_196             & 0.747425                   & 0.768764                   & 0.772191                   & 0.789858      & 0.798562      & 0.710980       & 0.610615         \\ \hline
auto\_price\_195         & 0.863399                   & 0.815007                   & 0.672240                   & 0.509488      & 0.747078      & 0.738998       & 0.693806         \\ \hline
bank32nh\_558            & 0.401371                   & 0.399995                   & 0.438601                   & 0.435068      & 0.419061      & 0.370140       & 0.247573         \\ \hline
bodyfat\_560             & 0.936925                   & 0.938313                   & 0.849033                   & 0.840143      & 0.948397      & 0.907751       & 0.858232         \\ \hline
cpu\_small\_227          & 0.912807                   & 0.923339                   & 0.927361                   & 0.922691      & 0.536490      & 0.878285       & 0.859981         \\ \hline
elevators\_216           & 0.410628                   & 0.452588                   & 0.494705                   & 0.499322      & 0.551810      & 0.387546       & 0.355308         \\ \hline
house\_16H\_574          & 0.365959                   & 0.400449                   & 0.404487                   & 0.400990      & 0.334073      & 0.339442       & 0.248186         \\ \hline
houses\_537              & 0.534771                   & 0.548320                   & 0.586357                   & 0.584077      & 0.563709      & 0.494428       & 0.472886         \\ \hline
kin8nm\_189              & 0.389164                   & 0.405103                   & 0.460511                   & 0.456236      & 0.393324      & 0.340757       & 0.294860         \\ \hline
mtp\_405                 & 0.318501                   & 0.336691                   & 0.357685                   & 0.354267      & 0.359294      & 0.296157       & 0.183265         \\ \hline
mv\_344                  & 0.914857                   & 0.928591                   & 0.930240                   & 0.937472      & 0.963350      & 0.847946       & 0.890945         \\ \hline
no2\_547                 & 0.428888                   & 0.479166                   & 0.472051                   & 0.458165      & 0.516820      & 0.381181       & 0.382920         \\ \hline
pol\_201                 & 0.750334                   & 0.723749                   & 0.779218                   & 0.772502      & 0.198512      & 0.700354       & 0.569024         \\ \hline
pollen\_529              & 0.521309                   & 0.486982                   & 0.601252                   & 0.591175      & 0.558306      & 0.438864       & 0.380181         \\ \hline
puma32H\_308             & 0.656846                   & 0.651822                   & 0.669697                   & 0.665608      & 0.204270      & 0.451785       & 0.606380         \\ \hline
satellite\_image\_294    & 0.706395                   & 0.748535                   & 0.744949                   & 0.743465      & 0.733783      & 0.670141       & 0.452368         \\ \hline
socmob\_541              & 0.669874                   & 0.611381                   & 0.637742                   & 0.608317      & 0.677353      & 0.582189       & 0.643831         \\ \hline
space\_ga\_507           & 0.502022                   & 0.520441                   & 0.454429                   & 0.445174      & 0.520454      & 0.439996       & 0.405508         \\ \hline
stock\_223               & 0.885754                   & 0.907392                   & 0.889649                   & 0.905912      & 0.850638      & 0.818754       & 0.865049         \\ \hline
tecator\_505             & 0.959962                   & 0.954073                   & 0.827610                   & 0.850458      & 0.972255      & 0.946702       & 0.857634         \\ \hline
us\_crime\_315           & 0.578675                   & 0.597916                   & 0.607469                   & 0.602808      & 0.351041      & 0.553025       & 0.372902         \\ \hline
vinnie\_519              & 0.698352                   & 0.704060                   & 0.608716                   & 0.635216      & 0.716773      & 0.596746       & 0.706372         \\ \hline
visualizing\_galaxy\_690 & 0.930156                   & 0.947740                   & 0.837415                   & 0.871342      & 0.957811      & 0.876663       & 0.913459         \\ \hline
wind\_503                & 0.669185                   & 0.670398                   & 0.685995                   & 0.679164      & 0.746188      & 0.642927       & 0.410033         \\ \hline
wine\_quality\_287       & 0.278516                   & 0.272298                   & 0.307583                   & 0.309223      & 0.291475      & 0.246417       & 0.227999         \\ \hline
\rowcolor[HTML]{FFFFC7} 
\textbf{Average Rank}    & \textbf{3.9}                        & \textbf{3.4}                        & \textbf{2.8}                        & \textbf{3.1}           & \textbf{3.0}           & \textbf{5.7}            & \textbf{6.0}              \\ \hline
\end{tabular}
\vspace{5mm}
\caption{Average Test $R^2$ by dataset: We show here the out-of-sample performance, measured using test $R^2$, for all datasets in our experiment. These results are the averages, obtained over a 10-fold CV. Standard errors are shown in Table \ref{performance_stds.table}}. \label{performance_means.table}
\end{table*}

\begin{table*}[h]
\begin{tabular}{|c|r|r|r|r|r|r|r|}
\hline
\textbf{Dataset}         & \multicolumn{1}{c|}{\textbf{MOSS-$\epsilon$-H}} & \multicolumn{1}{c|}{\textbf{MOSS-$\epsilon$-M}} & \multicolumn{1}{c|}{\textbf{MOSS-$\epsilon$-L}} & \multicolumn{1}{c|}{\textbf{FIRE}} & \multicolumn{1}{c|}{\textbf{GLRM}} & \multicolumn{1}{c|}{\textbf{SIRUS}} & \multicolumn{1}{c|}{\textbf{RULEFIT}} \\ \hline
Ailerons\_296            & 0.002285                                        & 0.002516                                        & 0.002836                                        & 0.002821                           & 0.00415                            & 0.002121                            & 0.00287                               \\ \hline
ESL\_1027                & 0.005362                                        & 0.004986                                        & 0.004822                                        & 0.003459                           & 0.003859                           & 0.00531                             & 0.009543                              \\ \hline
mercedes\_42570          & 0.007692                                        & 0.006812                                        & 0.00813                                         & 0.006828                           & 0.000123                           & 0.007624                            & 0.007663                              \\ \hline
Moneyball\_41021         & 0.001384                                        & 0.001677                                        & 0.001983                                        & 0.001883                           & 0.000951                           & 0.001561                            & 0.005223                              \\ \hline
abalone\_183             & 0.004065                                        & 0.003866                                        & 0.002029                                        & 0.00233                            & 0.001936                           & 0.00325                             & 0.003582                              \\ \hline
autoMpg\_196             & 0.007127                                        & 0.007724                                        & 0.004803                                        & 0.007234                           & 0.007407                           & 0.006754                            & 0.005385                              \\ \hline
auto\_price\_195         & 0.004301                                        & 0.012749                                        & 0.020409                                        & 0.023252                           & 0.011682                           & 0.014493                            & 0.013724                              \\ \hline
bank32nh\_558            & 0.002263                                        & 0.002518                                        & 0.001691                                        & 0.001815                           & 0.0034                             & 0.002477                            & 0.00371                               \\ \hline
bodyfat\_560             & 0.003717                                        & 0.004767                                        & 0.008932                                        & 0.009656                           & 0.002343                           & 0.002949                            & 0.005828                              \\ \hline
cpu\_small\_227          & 0.001989                                        & 0.001869                                        & 0.001205                                        & 0.001792                           & 0.002345                           & 0.002134                            & 0.003018                              \\ \hline
elevators\_216           & 0.002923                                        & 0.002512                                        & 0.002744                                        & 0.001783                           & 0.006588                           & 0.002986                            & 0.002997                              \\ \hline
house\_16H\_574          & 0.002709                                        & 0.002277                                        & 0.003425                                        & 0.003237                           & 0.001564                           & 0.002798                            & 0.003822                              \\ \hline
houses\_537              & 0.001252                                        & 0.00136                                         & 0.002994                                        & 0.002617                           & 0.002434                           & 0.001083                            & 0.001938                              \\ \hline
kin8nm\_189              & 0.002086                                        & 0.001958                                        & 0.002858                                        & 0.00196                            & 0.00322                            & 0.002149                            & 0.000992                              \\ \hline
mtp\_405                 & 0.005264                                        & 0.005478                                        & 0.003156                                        & 0.003498                           & 0.003801                           & 0.004378                            & 0.002034                              \\ \hline
mv\_344                  & 0.000899                                        & 0.000281                                        & 0.000253                                        & 0.000357                           & 0.000604                           & 0.00063                             & 0.001155                              \\ \hline
no2\_547                 & 0.014793                                        & 0.012096                                        & 0.011815                                        & 0.013442                           & 0.010413                           & 0.017488                            & 0.010606                              \\ \hline
pol\_201                 & 0.001823                                        & 0.003061                                        & 0.002274                                        & 0.001908                           & 0.002216                           & 0.003823                            & 0.003662                              \\ \hline
pollen\_529              & 0.002837                                        & 0.005641                                        & 0.004749                                        & 0.004218                           & 0.004445                           & 0.005575                            & 0.00268                               \\ \hline
puma32H\_308             & 0.001736                                        & 0.002041                                        & 0.001943                                        & 0.001705                           & 0.002937                           & 0.009204                            & 0.005404                              \\ \hline
satellite\_image\_294    & 0.001982                                        & 0.00186                                         & 0.00152                                         & 0.001805                           & 0.004031                           & 0.002496                            & 0.006522                              \\ \hline
socmob\_541              & 0.009706                                        & 0.009399                                        & 0.02505                                         & 0.021329                           & 0.007091                           & 0.01136                             & 0.005021                              \\ \hline
space\_ga\_507           & 0.006704                                        & 0.006642                                        & 0.020675                                        & 0.020141                           & 0.00524                            & 0.007398                            & 0.006961                              \\ \hline
stock\_223               & 0.003355                                        & 0.002338                                        & 0.002685                                        & 0.001824                           & 0.004715                           & 0.004208                            & 0.002727                              \\ \hline
tecator\_505             & 0.001079                                        & 0.001199                                        & 0.0135                                          & 0.010467                           & 0.000717                           & 0.001254                            & 0.00656                               \\ \hline
us\_crime\_315           & 0.005858                                        & 0.006185                                        & 0.00602                                         & 0.005627                           & 0.006231                           & 0.007665                            & 0.005101                              \\ \hline
vinnie\_519              & 0.01316                                         & 0.011878                                        & 0.010362                                        & 0.011215                           & 0.009029                           & 0.018975                            & 0.011436                              \\ \hline
visualizing\_galaxy\_690 & 0.002371                                        & 0.002072                                        & 0.008285                                        & 0.006323                           & 0.002107                           & 0.004603                            & 0.002173                              \\ \hline
wind\_503                & 0.002286                                        & 0.002195                                        & 0.002025                                        & 0.002219                           & 0.001447                           & 0.002937                            & 0.004641                              \\ \hline
wine\_quality\_287       & 0.002742                                        & 0.001916                                        & 0.002367                                        & 0.002378                           & 0.00202                            & 0.002405                            & 0.002517                              \\ \hline
\end{tabular}
\vspace{5mm}
\caption{Standard Error of Test $R^2$ by dataset. These standard errors are obtained over a 10-fold CV. The mean test $R^2$ of our performance results are shown in Table \ref{performance_means.table}} \label{performance_stds.table}
\end{table*}

\begin{table*}[h]
\begin{tabular}{|c|c|c|c|c|c|c|c|}
\hline
\textbf{Dataset}         & \textbf{MOSS-$\epsilon$-H} & \textbf{MOSS-$\epsilon$-M} & \textbf{MOSS-$\epsilon$-L} & \textbf{FIRE} & \textbf{GLRM} & \textbf{SIRUS} & \textbf{RULEFIT} \\ \hline
Ailerons\_296            & 0.519373                   & 0.216654                   & 0.252604                   & 0.242296      & 0.342012      & 0.711917       & 0.230505         \\ \hline
ESL\_1027                & 0.539919                   & 0.422515                   & 0.505706                   & 0.232644      & 0.277109      & 0.655401       & 0.424257         \\ \hline
mercedes\_42570          & 0.622377                   & 0.493946                   & 0.426768                   & 0.317824      & 0.314512      & 0.673675       & 0.263956         \\ \hline
Moneyball\_41021         & 0.575336                   & 0.341164                   & 0.496633                   & 0.165632      & 0.245471      & 0.667643       & 0.256258         \\ \hline
abalone\_183             & 0.414021                   & 0.34291                    & 0.296429                   & 0.263169      & 0.365313      & 0.497989       & 0.18037          \\ \hline
autoMpg\_196             & 0.538893                   & 0.260317                   & 0.417147                   & 0.22853       & 0.145074      & 0.577619       & 0.161384         \\ \hline
auto\_price\_195         & 0.289273                   & 0.153148                   & 0.358583                   & 0.025644      & 0.05346       & 0.264154       & 0.076101         \\ \hline
bank32nh\_558            & 0.603915                   & 0.373862                   & 0.11825                    & 0.235648      & 0.622924      & 0.68127        & 0.411111         \\ \hline
bodyfat\_560             & 0.657985                   & 0.316479                   & 0.657835                   & 0.198854      & 0.188848      & 0.708759       & 0.193073         \\ \hline
cpu\_small\_227          & 0.632039                   & 0.424135                   & 0.366626                   & 0.323713      & 0.468433      & 0.798771       & 0.583913         \\ \hline
elevators\_216           & 0.459524                   & 0.399735                   & 0.200121                   & 0.329536      & 0.306148      & 0.591005       & 0.384127         \\ \hline
house\_16H\_574          & 0.680887                   & 0.324827                   & 0.299862                   & 0.271005      & 0.610905      & 0.777729       & 0.642267         \\ \hline
houses\_537              & 0.730639                   & 0.501701                   & 0.408506                   & 0.580154      & 0.657892      & 0.808856       & 0.665739         \\ \hline
kin8nm\_189              & 0.34328                    & 0.39418                    & 0.174172                   & 0.2484        & 0.189194      & 0.376931       & 0.364127         \\ \hline
mtp\_405                 & 0.544872                   & 0.265218                   & 0.151811                   & 0.168422      & 0.109098      & 0.629402       & 0.332087         \\ \hline
mv\_344                  & 0.682011                   & 0.631429                   & 0.766772                   & 0.546314      & 0.650213      & 0.878519       & 0.693439         \\ \hline
no2\_547                 & 0.292112                   & 0.283675                   & 0.2304                     & 0.170359      & 0.24391       & 0.295515       & 0.281209         \\ \hline
pol\_201                 & 0.756046                   & 0.567302                   & 0.60816                    & 0.698832      & 0.518501      & 0.814867       & 0.646029         \\ \hline
pollen\_529              & 0.685926                   & 0.587831                   & 0.331494                   & 0.398385      & 0.413704      & 0.686772       & 0.477884         \\ \hline
puma32H\_308             & 0.632764                   & 0.810834                   & 0.405413                   & 0.386241      & 0.119466      & 0.664029       & 0.679414         \\ \hline
satellite\_image\_294    & 0.605291                   & 0.334127                   & 0.331583                   & 0.229861      & 0.347707      & 0.691534       & 0.313161         \\ \hline
socmob\_541              & 0.695397                   & 0.654392                   & 0.630077                   & 0.557368      & 0.170435      & 0.651693       & 0.533333         \\ \hline
space\_ga\_507           & 0.485637                   & 0.538771                   & 0.204127                   & 0.349709      & 0.581781      & 0.514082       & 0.363561         \\ \hline
stock\_223               & 0.579577                   & 0.541905                   & 0.471184                   & 0.343286      & 0.618864      & 0.697354       & 0.430476         \\ \hline
tecator\_505             & 0.656032                   & 0.253386                   & 0.805942                   & 0.079609      & 0.073071      & 0.746772       & 0.159471         \\ \hline
us\_crime\_315           & 0.428726                   & 0.234546                   & 0.215438                   & 0.034593      & 0.131451      & 0.518055       & 0.14915          \\ \hline
vinnie\_519              & 0.641058                   & 0.566561                   & 0.62381                    & 0.293618      & 0.47808       & 0.699577       & 0.346772         \\ \hline
visualizing\_galaxy\_690 & 0.594709                   & 0.424127                   & 0.526939                   & 0.182653      & 0.106043      & 0.713228       & 0.288042         \\ \hline
wind\_503                & 0.524868                   & 0.237302                   & 0.326007                   & 0.192867      & 0.300686      & 0.611111       & 0.311508         \\ \hline
wine\_quality\_287       & 0.469499                   & 0.443488                   & 0.247871                   & 0.345134      & 0.390176      & 0.530989       & 0.325507         \\ \hline
\rowcolor[HTML]{FFFFC7} 
\textbf{Average Rank}    & \textbf{2.3}               & \textbf{3.9}               & \textbf{4.6}               & \textbf{5.9}  & \textbf{5.3}  & \textbf{1.3}   & \textbf{4.6}     \\ \hline
\end{tabular} 
\vspace{5mm}
\caption{Average Empirical Stability (Measured using Average Pairwise DSC) by dataset: We show here the average empirical stability for all datasets in our experiments. These estimates were obtained over a 10-fold CV. Standard errors are shown in Table \ref{stability_std.table}}\label{stability_means.table}
\end{table*}

\begin{table*}[]
\begin{tabular}{|c|r|r|r|r|r|r|r|}
\hline
\textbf{Dataset}         & \multicolumn{1}{c|}{\textbf{MOSS-$\epsilon$-H}} & \multicolumn{1}{c|}{\textbf{MOSS-$\epsilon$-M}} & \multicolumn{1}{c|}{\textbf{MOSS-$\epsilon$-L}} & \multicolumn{1}{c|}{\textbf{FIRE}} & \multicolumn{1}{c|}{\textbf{GLRM}} & \multicolumn{1}{c|}{\textbf{SIRUS}} & \multicolumn{1}{c|}{\textbf{RULEFIT}} \\ \hline
Ailerons\_296            & 0.118842                                        & 0.079678                                        & 0.102419                                        & 0.104818                           & 0.123415                           & 0.120866                            & 0.088508                              \\ \hline
ESL\_1027                & 0.172018                                        & 0.164152                                        & 0.117612                                        & 0.111256                           & 0.20201                            & 0.238199                            & 0.169534                              \\ \hline
mercedes\_42570          & 0.084819                                        & 0.078256                                        & 0.084641                                        & 0.137956                           & 0.13561                            & 0.13074                             & 0.083186                              \\ \hline
Moneyball\_41021         & 0.109747                                        & 0.098393                                        & 0.148511                                        & 0.094008                           & 0.208727                           & 0.094479                            & 0.10799                               \\ \hline
abalone\_183             & 0.132611                                        & 0.118806                                        & 0.12475                                         & 0.114168                           & 0.177402                           & 0.149801                            & 0.090511                              \\ \hline
autoMpg\_196             & 0.123852                                        & 0.074444                                        & 0.180301                                        & 0.132714                           & 0.100081                           & 0.151147                            & 0.109536                              \\ \hline
auto\_price\_195         & 0.159446                                        & 0.091233                                        & 0.138117                                        & 0.040316                           & 0.093504                           & 0.16397                             & 0.081988                              \\ \hline
bank32nh\_558            & 0.090781                                        & 0.084361                                        & 0.116622                                        & 0.171948                           & 0.091741                           & 0.114701                            & 0.104743                              \\ \hline
bodyfat\_560             & 0.100415                                        & 0.091985                                        & 0.093507                                        & 0.10126                            & 0.094722                           & 0.113837                            & 0.112233                              \\ \hline
cpu\_small\_227          & 0.107813                                        & 0.085746                                        & 0.143522                                        & 0.119897                           & 0.223344                           & 0.098827                            & 0.080084                              \\ \hline
elevators\_216           & 0.181387                                        & 0.112785                                        & 0.134456                                        & 0.142449                           & 0.137291                           & 0.219673                            & 0.101059                              \\ \hline
house\_16H\_574          & 0.065372                                        & 0.105668                                        & 0.101722                                        & 0.143274                           & 0.104602                           & 0.087089                            & 0.076819                              \\ \hline
houses\_537              & 0.068961                                        & 0.107299                                        & 0.130898                                        & 0.113962                           & 0.144677                           & 0.070389                            & 0.109541                              \\ \hline
kin8nm\_189              & 0.135261                                        & 0.08342                                         & 0.09613                                         & 0.119301                           & 0.107888                           & 0.170825                            & 0.117471                              \\ \hline
mtp\_405                 & 0.159357                                        & 0.095143                                        & 0.099004                                        & 0.104523                           & 0.095361                           & 0.172323                            & 0.133667                              \\ \hline
mv\_344                  & 0.073298                                        & 0.099027                                        & 0.119891                                        & 0.11089                            & 0.167417                           & 0.061751                            & 0.073664                              \\ \hline
no2\_547                 & 0.14432                                         & 0.092084                                        & 0.134929                                        & 0.110993                           & 0.137453                           & 0.144827                            & 0.112705                              \\ \hline
pol\_201                 & 0.06828                                         & 0.078275                                        & 0.212616                                        & 0.107391                           & 0.265527                           & 0.082779                            & 0.111208                              \\ \hline
pollen\_529              & 0.078259                                        & 0.05903                                         & 0.181264                                        & 0.101756                           & 0.170623                           & 0.144676                            & 0.101434                              \\ \hline
puma32H\_308             & 0.127927                                        & 0.081662                                        & 0.108145                                        & 0.143083                           & 0.061426                           & 0.120889                            & 0.110327                              \\ \hline
satellite\_image\_294    & 0.106714                                        & 0.109409                                        & 0.154514                                        & 0.123443                           & 0.15882                            & 0.113462                            & 0.129678                              \\ \hline
socmob\_541              & 0.097516                                        & 0.094177                                        & 0.12578                                         & 0.10664                            & 0.066391                           & 0.102278                            & 0.107464                              \\ \hline
space\_ga\_507           & 0.112121                                        & 0.096492                                        & 0.137084                                        & 0.110447                           & 0.186371                           & 0.094293                            & 0.091464                              \\ \hline
stock\_223               & 0.090338                                        & 0.115773                                        & 0.105315                                        & 0.118745                           & 0.176057                           & 0.111832                            & 0.078487                              \\ \hline
tecator\_505             & 0.070615                                        & 0.071412                                        & 0.086415                                        & 0.07219                            & 0.084201                           & 0.11017                             & 0.082183                              \\ \hline
us\_crime\_315           & 0.188118                                        & 0.097784                                        & 0.105851                                        & 0.041837                           & 0.120541                           & 0.185771                            & 0.093401                              \\ \hline
vinnie\_519              & 0.103782                                        & 0.085893                                        & 0.126296                                        & 0.152109                           & 0.164701                           & 0.132472                            & 0.139373                              \\ \hline
visualizing\_galaxy\_690 & 0.124876                                        & 0.106514                                        & 0.139881                                        & 0.096474                           & 0.156731                           & 0.114932                            & 0.094652                              \\ \hline
wind\_503                & 0.120254                                        & 0.118262                                        & 0.107745                                        & 0.133936                           & 0.183295                           & 0.129051                            & 0.112851                              \\ \hline
wine\_quality\_287       & 0.105566                                        & 0.098378                                        & 0.127814                                        & 0.142321                           & 0.189997                           & 0.146662                            & 0.118966                              \\ \hline
\end{tabular}
\vspace{5mm}
\caption{Standard Error of Empirical Stability by dataset: We show here the standard error of our stability estimates, computed across a 10-fold CV. We report average stability in Table \ref{stability_means.table}}\label{stability_std.table}
\end{table*}

\subsection{Experiment Results: Comparisons} \label{comparisons.appx}
One 24 out of 30 of the datasets in our experiment, $MOSS$-$\epsilon$-H and SIRUS have empirical stabilities that are within 1 standard error of each other. These datasets are:

\begin{itemize}
 \item                                 ESL
 \item       mercedes
 \item                                  Moneyball
 \item                                     abalone
 \item                                     autoMpg
 \item                                   auto\_price
 \item                                    bank32nh
 \item                                     bodyfat
 \item                                  elevators
 \item                                      kin8nm
 \item                                        mtp
 \item                                         no2
 \item                                        pol
 \item                                      pollen
 \item                                     puma32H
 \item                             satellite\_image
 \item                                     socmob
 \item                                    space\_ga
 \item                                    tecator
 \item                                   us\_crime
 \item                                      vinnie
 \item                          visualizing\_galaxy
 \item                                        wind
 \item                                wine\_quality
\end{itemize}

\begin{figure}[h]
    \centering
    \includegraphics[width=.95\linewidth]{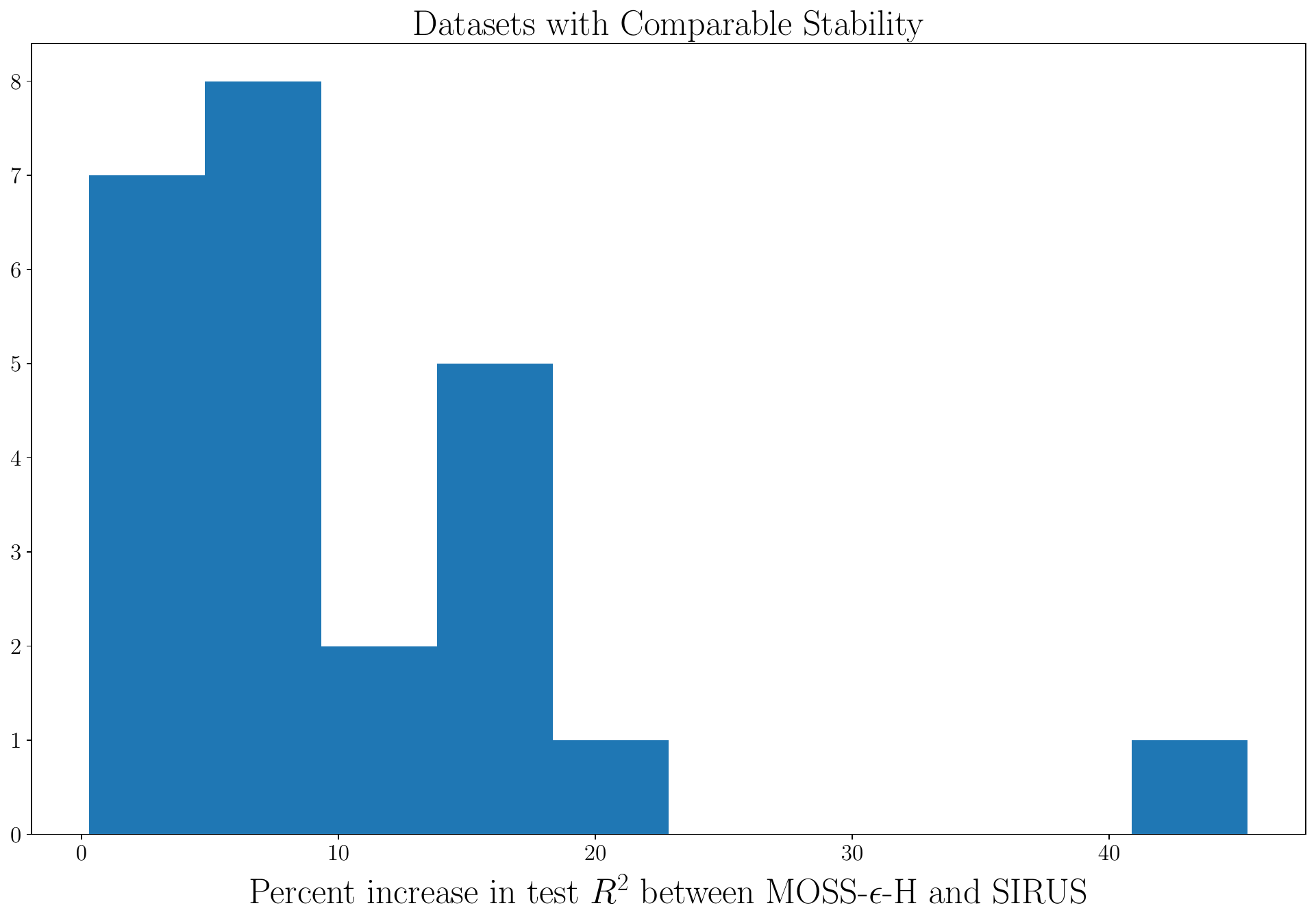}
    \caption{Percent increase in test $R^2$ between MOSS-$\epsilon$-H and SIRUS on datasets where the stability of the both methods is comparable.}
    \label{d1.fig}
\end{figure}

On these datasets, we observe than MOSS-$\epsilon$-H produces much more accurate rule sets, with an 10 percent increase of 10$R^2$. We show the distribution of this percent increase in test performance, between MOSS-$\epsilon$-H and SIRUS, for these datasets where the rule sets produced have comparable empirical stability, Figure \ref{d1.fig}.

\begin{figure}[h]
    \centering
    \includegraphics[width=.95\linewidth]{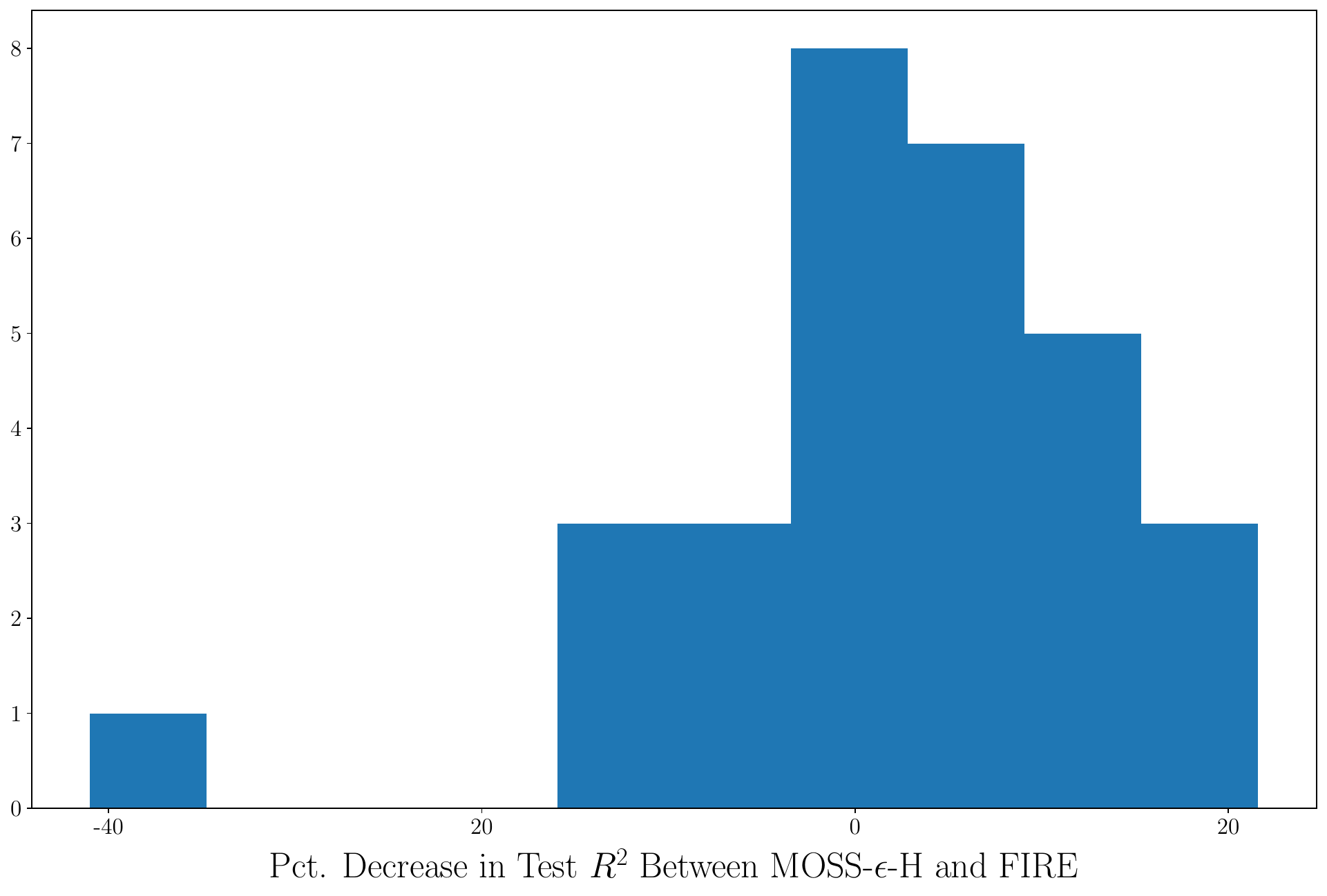}
    \caption{Percent decrease in test $R^2$ between MOSS-$\epsilon$-H and FIRE.}
    \label{d2.fig}
\end{figure}

We also compare MOSS-$\epsilon$-H against FIRE, the competing method with the best accuracy. We observe that the on average, $MOSS-\epsilon$-H has a 2 percent decrease in test $R^2$ compared to FIRE, we show the distribution in Figure \ref{d2.fig}.

\begin{figure}[h]
    \centering
    \includegraphics[width=.95\linewidth]{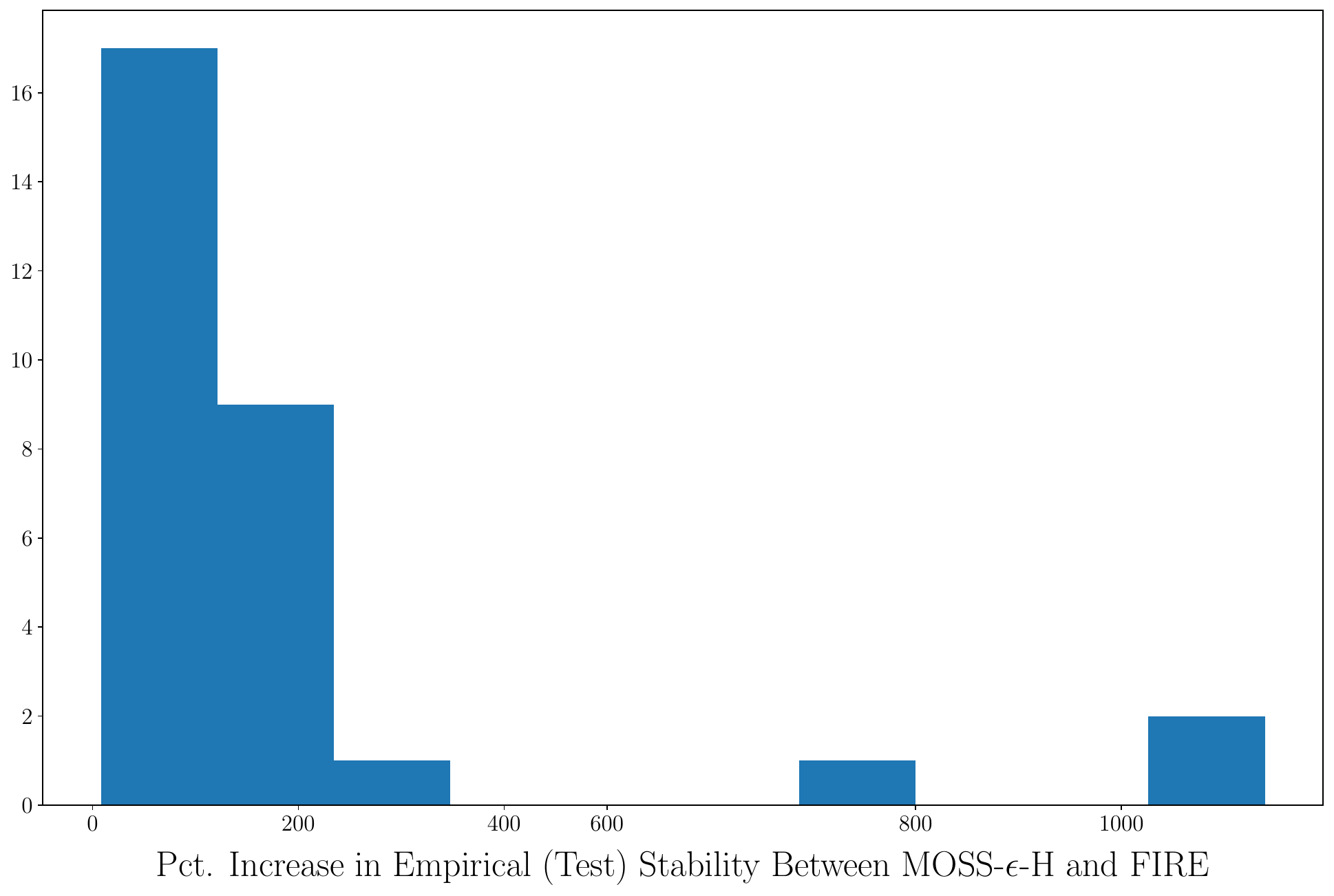}
    \caption{Percent increase in empirical (test) stability between MOSS-$\epsilon$-H and FIRE.}
    \label{d3.fig}
\end{figure}

However, the models produced by $MOSS$-$\epsilon$-H are much more stable, with an average percent increase in empirical stability (our average pairwise DSC metric) of 190 percent. We show the distribution of this if Figure \ref{d3.fig}.

\clearpage
\newpage

\subsection{Sensitivity Analyses} \label{sens.section}
In this section, we analyze the sensitivity of MOSS to parameters $\gamma$ and $k$.

\subsubsection{\textbf{Parameter $\gamma$}} \label{gammasense.section} This parameter controls the ridge regularization penalty in accuracy objective $H_2(z)$. We follow the procedure below to assess the sensitivity of MOSS to $\gamma$. 

We use the 30 OpenML datasets from our experiments in \S\ref{experiments.section} and repeat a 10-fold CV on each dataset. On each training fold, we use random forests to generate $m \sim 10^3$ candidate decision rules. We then apply MOSS to construct rule sets with 15 decision rules. We repeat this procedure across all folds and datasets in our experiment while varying $\gamma \in \{0.0001, 0.0005,0.001,0.005\}$. For each value of $\gamma$ we record the average predictive performance (out-of-sample $R^2$) and empirical stability (average pairwise Dice-Sorensen coefficient) of the rule sets. We report the results of this sensitivity analysis in Table \ref{gamma_sens.table}. From this table, we observe that the accuracy of MOSS (out-of-sample $R^2$) is relatively \emph{insensitive} to $\gamma$. Across all datasets, the rules sets constructed with varying values of $\gamma$ achieve similar $R^2$ scores. We also note that increasing $\gamma$ appears to decrease the empirical stability of the rule sets, by a slight degree.

Our results here suggest that the performance of MOSS in terms of empirical stability and accuracy is relatively insensitive to parameter $\gamma$. As an aside, we also note that $\gamma$ influences the computation time of our cutting plane algorithm. When $\gamma$ is large, Algorithm \ref{OA.alg} requires more iterations to converge. As such, we recommend setting $\gamma$ to a small value around $10^{-2}$ or $10^{-3}$ as a default.

\subsubsection{\textbf{Parameter $k$}} \label{ksens.section} The parameter $k$ controls the size of the rule sets constructed by MOSS. It is important to restrict $k$ to be small so that the rule sets remain interpretable; in fact, \citep{bertsimas2021sparse} and \citep{liu2023fire} restrict rule sets to contain $< 20$ rules to remain human readable. We use the procedure below to assess the sensivity of MOSS to $k$

We employ the same 30 OpenML datasets from our experiments in \S\ref{experiments.section} and perform 10-fold cross-validation on each dataset.  On each training fold we generate $m~10^3$ candidate rules and then we apply MOSS to construct rule sets. We set $\gamma = 0.001$ vary $k \in \{5, 10, 15, 20, 25\}$. For each value of $k$ we record the average out-of-sample $R^2$ score and the empirical stability of the rule sets, averaged across all folds for each dataset. We report the results of this sensitivity analysis in Table \ref{k_sens.table}. 

From Table \ref{k_sens.table}, we observe that as $k$ increases, both the average out-of-sample $R^2$ and the empirical stability of the rule sets improve. Very compact rule sets ($k = 5$) exhibit significantly worse accuracy and stability compared to larger rule sets ($k>10$). While larger rule sets perform better in terms of accuracy and stability, they are less interpretable. Based on these results, we recommend setting $k$ to 10, 15, or 20 rules for MOSS to achieve a balance between accuracy, stability, and interpretability.

\begin{figure}[h]
    \centering
    \includegraphics[width=0.425\textwidth]{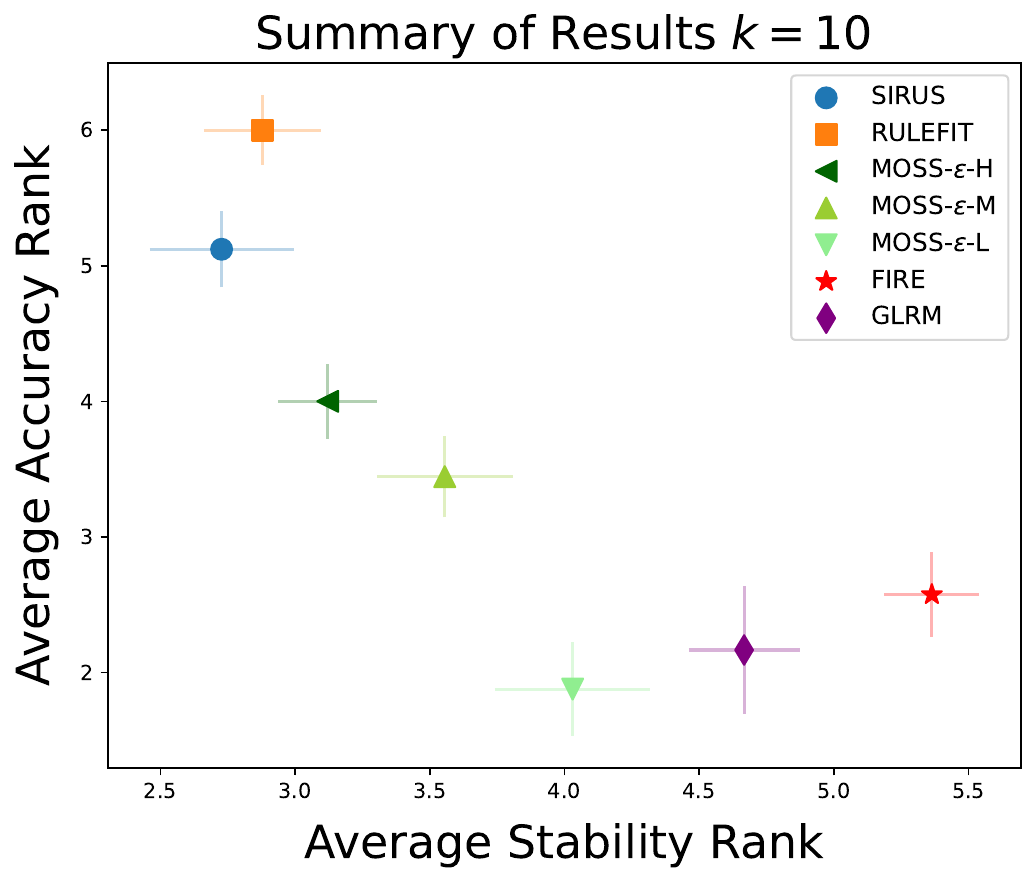}
    \caption{Experimental results for size $k = 10$ rule sets. MOSS methods are able to compute the Pareto frontier between accuracy and stability.}
    \label{k10.fig}
\end{figure}

\begin{figure}[h]
    \centering
    \includegraphics[width=0.425\textwidth]{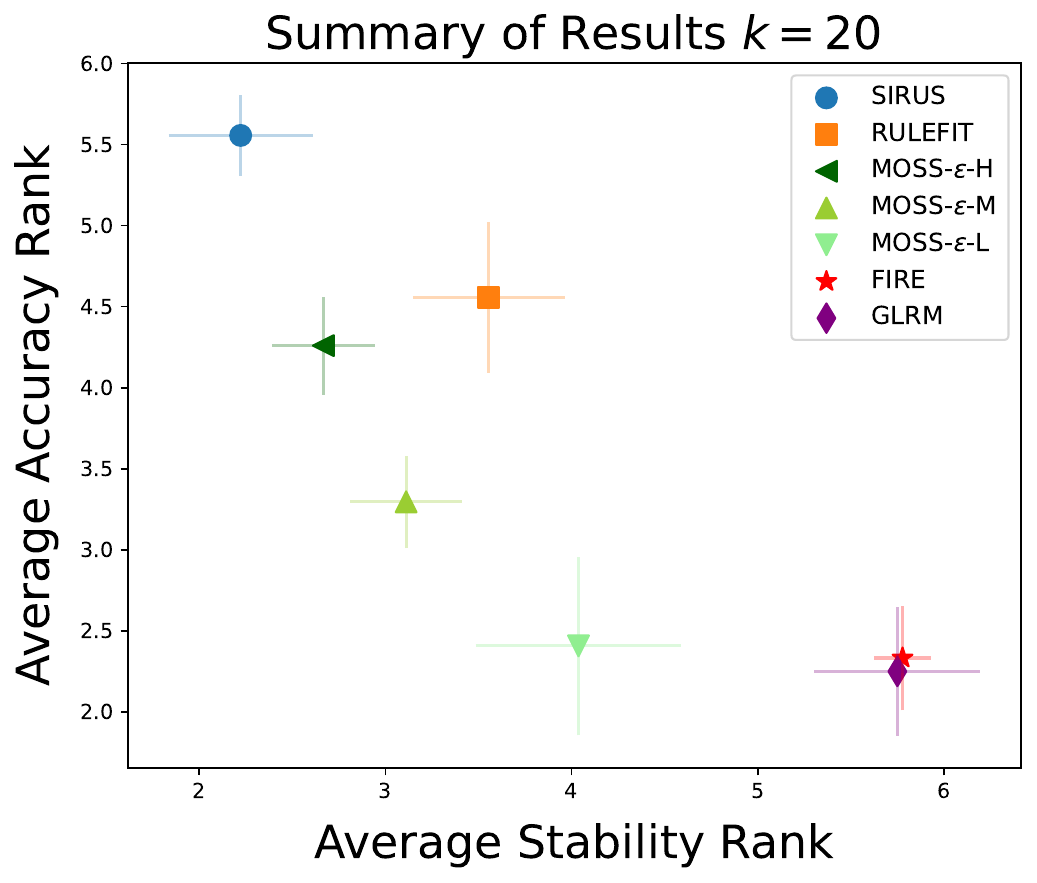}
    \caption{Experimental results for size $k = 20$ rule sets.}
    \label{k20.fig}
\end{figure}

Our current experimental results in \S\ref{experiments.section} compare MOSS against competing algorithms for $k = 15$ sized rule sets. We repeat our experiments for $k = 10$ and $k = 20$; Figures \ref{k10.fig} and \ref{k20.fig} show our results for these experiments. From these plots, we again see that MOSS outperforms our competing algorithms jointly in terms of both accuracy and stability, and achieves a balance between the two objectives.

\newpage
\twocolumn
\tablefirsthead{%
  \hline
  Dataset Name & $\gamma$ & $R^2$ & Stability \\
  \hline
}
\tablehead{%
  \hline
  Dataset Name & $\gamma$ & $R^2$ & Stability \\
  \hline
}
\tablelasttail{\hline}
\topcaption{Results for sensitivity analysis over parameter $\gamma$ in MOSS.}
\label{gamma_sens.table}
\begin{supertabular}{|l|c|c|c|}
auto\_price & 0.0001 & 0.6637 & 0.475\\
 & 0.0005 & 0.6348 & 0.4765\\
 & 0.001 & 0.6455 & 0.47\\
 & 0.005 & 0.6472 & 0.4763\\
tecator & 0.0001 & 0.8987 & 0.7186\\
 & 0.0005 & 0.8988 & 0.7217\\
 & 0.001 & 0.8988 & 0.7217\\
 & 0.005 & 0.8984 & 0.6922\\
body\_fat & 0.0001 & 0.8665 & 0.5647\\
 & 0.0005 & 0.8648 & 0.5612\\
 & 0.001 & 0.865 & 0.5559\\
 & 0.005 & 0.8678 & 0.5465\\
visualizing\_galaxy & 0.0001 & 0.8134 & 0.5652\\
 & 0.0005 & 0.8203 & 0.545\\
 & 0.001 & 0.8133 & 0.5275\\
 & 0.005 & 0.8295 & 0.4612\\
vinnie & 0.0001 & 0.5648 & 0.6211\\
 & 0.0005 & 0.5639 & 0.6509\\
 & 0.001 & 0.5641 & 0.632\\
 & 0.005 & 0.5648 & 0.6039\\
autoMpg & 0.0001 & 0.6703 & 0.5688\\
 & 0.0005 & 0.6703 & 0.5654\\
 & 0.001 & 0.6695 & 0.5638\\
 & 0.005 & 0.6759 & 0.5366\\
ESL & 0.0001 & 0.7574 & 0.4205\\
 & 0.0005 & 0.7566 & 0.4285\\
 & 0.001 & 0.7611 & 0.4121\\
 & 0.005 & 0.7731 & 0.4076\\
no2 & 0.0001 & 0.3806 & 0.3293\\
 & 0.0005 & 0.3787 & 0.3054\\
 & 0.001 & 0.3815 & 0.3114\\
 & 0.005 & 0.3892 & 0.3158\\
stock & 0.0001 & 0.8213 & 0.7335\\
 & 0.0005 & 0.8284 & 0.7116\\
 & 0.001 & 0.8367 & 0.7262\\
 & 0.005 & 0.8414 & 0.7045\\
socmob & 0.0001 & 0.5895 & 0.7789\\
 & 0.0005 & 0.5837 & 0.7759\\
 & 0.001 & 0.5951 & 0.7654\\
 & 0.005 & 0.613 & 0.7337\\
Moneyball & 0.0001 & 0.8333 & 0.4698\\
 & 0.0005 & 0.8324 & 0.4547\\
 & 0.001 & 0.8346 & 0.4707\\
 & 0.005 & 0.8419 & 0.4088\\
us\_crime & 0.0001 & 0.5371 & 0.561\\
 & 0.0005 & 0.5607 & 0.5553\\
 & 0.001 & 0.5603 & 0.5191\\
 & 0.005 & 0.5597 & 0.5388\\
 \end{supertabular}
 \pagebreak
 \tablelasttail{\hline}
 \begin{supertabular}{|l|c|c|c|}
space\_ga & 0.0001 & 0.4742 & 0.4241\\
 & 0.0005 & 0.4776 & 0.4316\\
 & 0.001 & 0.4781 & 0.4264\\
 & 0.005 & 0.4817 & 0.384\\
pollen & 0.0001 & 0.4541 & 0.6702\\
 & 0.0005 & 0.4632 & 0.6729\\
 & 0.001 & 0.4627 & 0.6654\\
 & 0.005 & 0.4784 & 0.6423\\
abalone & 0.0001 & 0.3636 & 0.3684\\
 & 0.0005 & 0.3802 & 0.3317\\
 & 0.001 & 0.3824 & 0.3484\\
 & 0.005 & 0.3777 & 0.3434\\
Mercedes\_Benz & 0.0001 & 0.564 & 0.737\\
 & 0.0005 & 0.5724 & 0.6888\\
 & 0.001 & 0.5728 & 0.6523\\
 & 0.005 & 0.5735 & 0.6444\\
mtp & 0.0001 & 0.2939 & 0.5516\\
 & 0.0005 & 0.2966 & 0.5172\\
 & 0.001 & 0.2979 & 0.5126\\
 & 0.005 & 0.3022 & 0.4279\\
satellite\_image & 0.0001 & 0.6813 & 0.6698\\
 & 0.0005 & 0.6836 & 0.6254\\
 & 0.001 & 0.6992 & 0.6063\\
 & 0.005 & 0.7006 & 0.6397\\
wine\_quality & 0.0001 & 0.2578 & 0.5483\\
 & 0.0005 & 0.2588 & 0.537\\
 & 0.001 & 0.2586 & 0.5281\\
 & 0.005 & 0.2564 & 0.5564\\
wind & 0.0001 & 0.6522 & 0.5411\\
 & 0.0005 & 0.6551 & 0.5262\\
 & 0.001 & 0.6598 & 0.5325\\
 & 0.005 & 0.6601 & 0.5139\\
puma32H & 0.0001 & 0.4491 & 0.5895\\
 & 0.0005 & 0.4693 & 0.5849\\
 & 0.001 & 0.4668 & 0.6016\\
 & 0.005 & 0.4943 & 0.6027\\
bank32nh & 0.0001 & 0.3444 & 0.6355\\
 & 0.0005 & 0.361 & 0.59\\
 & 0.001 & 0.3618 & 0.5591\\
 & 0.005 & 0.391 & 0.6066\\
cpu\_small & 0.0001 & 0.8902 & 0.7885\\
 & 0.0005 & 0.8966 & 0.7579\\
 & 0.001 & 0.9022 & 0.7471\\
 & 0.005 & 0.9057 & 0.744\\
kin8nm & 0.0001 & 0.3544 & 0.4219\\
 & 0.0005 & 0.3447 & 0.3878\\
 & 0.001 & 0.3602 & 0.3835\\
 & 0.005 & 0.3801 & 0.3585\\
 \end{supertabular}
\tablelasttail{\hline}
\begin{supertabular}{|l|c|c|c|}
Ailerons & 0.0001 & 0.5956 & 0.6722\\
 & 0.0005 & 0.606 & 0.6721\\
 & 0.001 & 0.6114 & 0.6802\\
 & 0.005 & 0.6119 & 0.6747\\
pol & 0.0001 & 0.6868 & 0.8923\\
 & 0.0005 & 0.7077 & 0.8658\\
 & 0.001 & 0.7345 & 0.818\\
 & 0.005 & 0.7385 & 0.8203\\
elevators & 0.0001 & 0.4042 & 0.5044\\
 & 0.0005 & 0.4056 & 0.4557\\
 & 0.001 & 0.4085 & 0.4521\\
 & 0.005 & 0.4185 & 0.4524\\
houses & 0.0001 & 0.5025 & 0.7706\\
 & 0.0005 & 0.5281 & 0.7889\\
 & 0.001 & 0.5304 & 0.7643\\
 & 0.005 & 0.5313 & 0.7493\\
house\_16H & 0.0001 & 0.3652 & 0.7906\\
 & 0.0005 & 0.3683 & 0.7537\\
 & 0.001 & 0.3674 & 0.7137\\
 & 0.005 & 0.3669 & 0.677\\
mv & 0.0001 & 0.8932 & 0.8677\\
 & 0.0005 & 0.9156 & 0.8006\\
 & 0.001 & 0.9171 & 0.8129\\
 & 0.005 & 0.9178 & 0.8049\\
\end{supertabular}

\clearpage
\newpage

\newpage
\twocolumn
\tablefirsthead{%
  \hline
  Dataset Name & $k$ & $R^2$ & Stability \\
  \hline
}
\tablehead{%
  \hline
  Dataset Name & $k$ & $R^2$ & Stability \\
  \hline
}
\tablelasttail{\hline}
\topcaption{Results for sensitivity analysis over parameter $k$ in MOSS.}
\label{k_sens.table}
\begin{supertabular}{|l|c|c|c|}
auto\_price  & 5 & 0.305 & 0.4418\\
 & 10 & 0.6095 & 0.5088\\
 & 15 & 0.6637 & 0.475\\
 & 20 & 0.6379 & 0.5171\\
 & 25 & 0.665 & 0.5274\\
tecator & 5 & 0.697 & 0.5608\\
 & 10 & 0.8939 & 0.7066\\
 & 15 & 0.8987 & 0.7186\\
 & 20 & 0.9036 & 0.6572\\
 & 25 & 0.9075 & 0.6601\\
bodyfat & 5 & 0.809 & 0.6278\\
 & 10 & 0.8583 & 0.6048\\
 & 15 & 0.8665 & 0.5647\\
 & 20 & 0.8655 & 0.5126\\
 & 25 & 0.8674 & 0.52\\
visualizing\_galaxy  & 5 & 0.6876 & 0.3067\\
 & 10 & 0.8018 & 0.5083\\
 & 15 & 0.8134 & 0.5652\\
 & 20 & 0.8447 & 0.5144\\
 & 25 & 0.8344 & 0.4785\\
vinnie  & 5 & 0.4535 & 0.2411\\
 & 10 & 0.5683 & 0.4605\\
 & 15 & 0.5648 & 0.6211\\
 & 20 & 0.5766 & 0.6514\\
 & 25 & 0.593 & 0.7375\\
autoMpg  & 5 & 0.585 & 0.3689\\
 & 10 & 0.6649 & 0.5461\\
 & 15 & 0.677 & 0.5721\\
 & 20 & 0.667 & 0.5395\\
 & 25 & 0.6647 & 0.5604\\
ESL & 5 & 0.6378 & 0.2533\\
 & 10 & 0.7222 & 0.3706\\
 & 15 & 0.7574 & 0.4205\\
 & 20 & 0.7609 & 0.3887\\
 & 25 & 0.7832 & 0.4377\\
no2  & 5 & 0.275 & 0.212\\
 & 10 & 0.3549 & 0.2908\\
 & 15 & 0.3806 & 0.3293\\
 & 20 & 0.3563 & 0.3614\\
 & 25 & 0.3625 & 0.3922\\
stock  & 5 & 0.7265 & 0.668\\
 & 10 & 0.8007 & 0.6531\\
 & 15 & 0.8213 & 0.7335\\
 & 20 & 0.8445 & 0.7309\\
 & 25 & 0.8417 & 0.715\\
socmob & 5 & 0.452 & 0.4278\\
 & 10 & 0.5742 & 0.6578\\
 & 15 & 0.5895 & 0.7789\\
 & 20 & 0.6511 & 0.8299\\
 & 25 & 0.6516 & 0.8772\\
 \end{supertabular}
 \pagebreak
 \tablelasttail{\hline}
 \begin{supertabular}{|l|c|c|c|}
Moneyball & 5 & 0.6956 & 0.3833\\
 & 10 & 0.8019 & 0.501\\
 & 15 & 0.8314 & 0.4918\\
 & 20 & 0.8488 & 0.5511\\
 & 25 & 0.8535 & 0.5291\\
us\_crime & 5 & 0.4951 & 0.4\\
 & 10 & 0.5301 & 0.5551\\
 & 15 & 0.5371 & 0.561\\
 & 20 & 0.5603 & 0.5638\\
 & 25 & 0.5646 & 0.5597\\
space\_ga  & 5 & 0.3851 & 0.3594\\
 & 10 & 0.4535 & 0.4059\\
 & 15 & 0.4742 & 0.4241\\
 & 20 & 0.475 & 0.4395\\
 & 25 & 0.4814 & 0.4664\\
pollen & 5 & 0.3377 & 0.5758\\
 & 10 & 0.3896 & 0.6755\\
 & 15 & 0.4541 & 0.6702\\
 & 20 & 0.4986 & 0.6835\\
 & 25 & 0.5001 & 0.6815\\
abalone & 5 & 0.2971 & 0.3056\\
 & 10 & 0.3518 & 0.3159\\
 & 15 & 0.3636 & 0.3684\\
 & 20 & 0.3725 & 0.4283\\
 & 25 & 0.3749 & 0.4625\\
Mercedes\_Benz & 5 & 0.4533 & 0.622\\
 & 10 & 0.5139 & 0.6763\\
 & 15 & 0.5596 & 0.7789\\
 & 20 & 0.5744 & 0.7735\\
 & 25 & 0.5744 & 0.7741\\
mtp & 5 & 0.2438 & 0.4286\\
 & 10 & 0.2849 & 0.4848\\
 & 15 & 0.2939 & 0.5516\\
 & 20 & 0.2967 & 0.5235\\
 & 25 & 0.305 & 0.5252\\
satellite\_image & 5 & 0.5849 & 0.413\\
 & 10 & 0.6647 & 0.6721\\
 & 15 & 0.6813 & 0.6698\\
 & 20 & 0.6927 & 0.6854\\
 & 25 & 0.6929 & 0.7176\\
wine\_quality & 5 & 0.223 & 0.5164\\
 & 10 & 0.2512 & 0.5585\\
 & 15 & 0.2578 & 0.5483\\
 & 20 & 0.2676 & 0.5849\\
 & 25 & 0.2745 & 0.6586\\
wind & 5 & 0.5662 & 0.4426\\
 & 10 & 0.6365 & 0.481\\
 & 15 & 0.6522 & 0.5411\\
 & 20 & 0.661 & 0.5683\\
 & 25 & 0.6674 & 0.6011\\
 \end{supertabular}
\tablelasttail{\hline}
\begin{supertabular}{|l|c|c|c|}
puma32H & 5 & 0.2995 & 0.5611\\
 & 10 & 0.3882 & 0.6184\\
 & 15 & 0.4491 & 0.5895\\
 & 20 & 0.4888 & 0.5856\\
 & 25 & 0.5198 & 0.5791\\
bank32nh & 5 & 0.2409 & 0.5289\\
 & 10 & 0.3364 & 0.617\\
 & 15 & 0.3444 & 0.6355\\
 & 20 & 0.3763 & 0.631\\
 & 25 & 0.3845 & 0.6459\\
cpu\_small & 5 & 0.7487 & 0.8\\
 & 10 & 0.8511 & 0.8267\\
 & 15 & 0.8902 & 0.7885\\
 & 20 & 0.8956 & 0.8143\\
 & 25 & 0.9017 & 0.8016\\
kin8nm & 5 & 0.2427 & 0.2961\\
 & 10 & 0.3249 & 0.377\\
 & 15 & 0.3544 & 0.4219\\
 & 20 & 0.3698 & 0.4399\\
 & 25 & 0.3826 & 0.4887\\
Ailerons & 5 & 0.5018 & 0.6537\\
 & 10 & 0.5707 & 0.6423\\
 & 15 & 0.5956 & 0.6722\\
 & 20 & 0.6105 & 0.7177\\
 & 25 & 0.6138 & 0.7105\\
pol & 5 & 0.5539 & 1.0\\
 & 10 & 0.668 & 0.9156\\
 & 15 & 0.6868 & 0.8923\\
 & 20 & 0.7157 & 0.9148\\
 & 25 & 0.747 & 0.9127\\
elevators & 5 & 0.3471 & 0.3426\\
 & 10 & 0.3945 & 0.5006\\
 & 15 & 0.4042 & 0.5044\\
 & 20 & 0.408 & 0.5216\\
 & 25 & 0.4114 & 0.5319\\
houses & 5 & 0.4159 & 0.7268\\
 & 10 & 0.4903 & 0.7339\\
 & 15 & 0.5025 & 0.7706\\
 & 20 & 0.5198 & 0.7719\\
 & 25 & 0.5336 & 0.7801\\
house\_16H & 5 & 0.2871 & 0.5215\\
 & 10 & 0.3413 & 0.7106\\
 & 15 & 0.3652 & 0.7906\\
 & 20 & 0.3708 & 0.7887\\
 & 25 & 0.3788 & 0.752\\
mv & 5 & 0.8005 & 0.816\\
 & 10 & 0.8467 & 0.8067\\
 & 15 & 0.9149 & 0.8391\\
 & 20 & 0.9214 & 0.8401\\
 & 25 & 0.9201 & 0.8424\\
\end{supertabular}

\clearpage
\newpage

\clearpage
\newpage

\subsection{Choice of Stability Measure} \label{choice_of_stability.section}

In \S\ref{model_stability_def.section} of the main text of our paper, we present our approach to asses the stability of rule algorithms. We use the Dice-Sorensens coefficient to measure the similarity between a pair of rule sets, $R_i$ and $R_j$, and we take the average pairwise coefficient across all rule sets to be the empirical stability of the algorithm.

We note here that we can use various alternative measures to assess the stability of our rule algorithms, and we show that our experimental results are \emph{insensitive} to our choice of stability measure. We evaluate the following metrics \cite{nogueira2018stability}.

\begin{align*}
    \text{Jaccard}(R_i, R_j) = \frac{|R_i \cap R_j|}{|R_i \cup R_j|}
\end{align*}

\begin{align*}
    \text{Ochiai}(R_i, R_j) = \frac{|R_i \cap R_j|}{\sqrt{|R_i|  |R_j|}}
\end{align*}

\begin{align*}
    \text{POG}(R_i, R_j) = \frac{|R_i \cap R_j|}{|R_i|}
\end{align*}

For each metric, we again take the average pairwise similarity between rule sets as the empirical stability of the rule algorithm. 

We repeat our experimental evaluation of MOSS (\S\ref{experiments.section}) and report our stability results using the measures discussed above. We summarize the results of this analysis in Table \ref{stab_measure.table}.

\begin{table}[h]
\scalebox{0.85}{
\begin{tabular}{|c|c|c|c|c|c|c|c|}
\hline
Metric  & MOSS-H & MOSS-M & MOSS-L & FIRE & GLRM & SIRUS & RuleFit \\ \hline
DSC     & 2.3    & 3.9    & 4.6    & 5.9  & 5.3  & 1.3   & 4.6     \\ \hline
Jaccard & 2.5    & 4.0    & 4.5    & 5.8  & 5.2  & 1.4   & 4.4     \\ \hline
Ochiai  & 2.5    & 3.8    & 4.6    & 5.8  & 5.3  & 1.4   & 4.3     \\ \hline
POG     & 2.4    & 3.9    & 4.6    & 5.9  & 5.2  & 1.4   & 4.3     \\ \hline
\end{tabular}}
\caption{Average stability ranking for each method in our experimental evaluation of MOSS using different stability measures.}
\label{stab_measure.table}
\end{table}

In this table, we present the average empirical stability ranking of each algorithm considered in our experiment, across all datasets. The first row reports stability measured using our Dice-Sorensen coefficient metrics, which represents the main result of our paper as shown in \S\ref{method_ranks.fig}. The subsequent rows display the stability rankings measured using our alternative metrics. This table demonstrates that our findings remain consistent regardless of the metrics used.

We report dataset-level stability results for each metric in the pages below (Tables \ref{jaccard.table},\ref{ochiai.table},\ref{pog.table}).

\clearpage
\newpage

\begin{table*}[h]
\begin{tabular}{|l|r|r|r|r|r|r|r|}
\hline
\multicolumn{1}{|c|}{\textbf{Dataset Name}}   & \multicolumn{1}{c|}{\textbf{MOSS-H}} & \multicolumn{1}{c|}{\textbf{MOSS-M}} & \multicolumn{1}{c|}{\textbf{MOSS-L}} & \multicolumn{1}{c|}{\textbf{FIRE}} & \multicolumn{1}{c|}{\textbf{GLRM}} & \multicolumn{1}{c|}{\textbf{SIRUS}} & \multicolumn{1}{c|}{\textbf{RuleFit}} \\ \hline
Ailerons\_296                                 & 0.37991                              & 0.129643                             & 0.159893                             & 0.142019                           & 0.13491                            & 0.601258                            & 0.13888                               \\ \hline
ESL\_1027                                     & 0.412596                             & 0.297479                             & 0.3629                               & 0.136341                           & 0.177267                           & 0.564619                            & 0.299636                              \\ \hline
Mercedes\_Benz\_42570 & 0.478343                             & 0.354673                             & 0.296266                             & 0.196763                           & 0.183194                           & 0.548374                            & 0.160499                              \\ \hline
Moneyball\_41021                              & 0.423393                             & 0.214538                             & 0.366631                             & 0.093294                           & 0.158887                           & 0.523579                            & 0.155207                              \\ \hline
abalone\_183                                  & 0.296459                             & 0.231639                             & 0.206711                             & 0.156669                           & 0.238788                           & 0.382027                            & 0.110963                              \\ \hline
autoMpg\_196                                  & 0.405238                             & 0.155535                             & 0.299937                             & 0.135751                           & 0.081399                           & 0.457141                            & 0.096783                              \\ \hline
auto\_price\_195                              & 0.197105                             & 0.09418                              & 0.251944                             & 0.013422                           & 0.030097                           & 0.178004                            & 0.044721                              \\ \hline
bank32nh\_558                                 & 0.451013                             & 0.238251                             & 0.071243                             & 0.14544                            & 0.459234                           & 0.543302                            & 0.270109                              \\ \hline
bodyfat\_560                                  & 0.536347                             & 0.202649                             & 0.54129                              & 0.114032                           & 0.107355                           & 0.604153                            & 0.117933                              \\ \hline
cpu\_small\_227                               & 0.495879                             & 0.285403                             & 0.244834                             & 0.199628                           & 0.332819                           & 0.716664                            & 0.432053                              \\ \hline
elevators\_216                                & 0.336396                             & 0.27418                              & 0.129023                             & 0.206686                           & 0.188657                           & 0.484291                            & 0.262566                              \\ \hline
house\_16H\_574                               & 0.547808                             & 0.207961                             & 0.187306                             & 0.164832                           & 0.448287                           & 0.682481                            & 0.50674                               \\ \hline
houses\_537                                   & 0.617404                             & 0.360401                             & 0.282309                             & 0.418224                           & 0.507125                           & 0.728404                            & 0.537911                              \\ \hline
kin8nm\_189                                   & 0.222827                             & 0.26181                              & 0.108051                             & 0.147727                           & 0.108419                           & 0.253536                            & 0.24064                               \\ \hline
mtp\_405                                      & 0.414122                             & 0.164794                             & 0.089602                             & 0.095773                           & 0.060559                           & 0.515025                            & 0.217215                              \\ \hline
mv\_344                                       & 0.542536                             & 0.479759                             & 0.657076                             & 0.384421                           & 0.505367                           & 0.82421                             & 0.557431                              \\ \hline
no2\_547                                      & 0.189344                             & 0.178384                             & 0.154659                             & 0.097312                           & 0.146761                           & 0.193289                            & 0.177325                              \\ \hline
pol\_201                                      & 0.686128                             & 0.442294                             & 0.508667                             & 0.5478                             & 0.396808                           & 0.782373                            & 0.538459                              \\ \hline
pollen\_529                                   & 0.540905                             & 0.428081                             & 0.224378                             & 0.253962                           & 0.275593                           & 0.553185                            & 0.325646                              \\ \hline
puma32H\_308                                  & 0.509926                             & 0.747302                             & 0.275336                             & 0.249453                           & 0.064664                           & 0.549834                            & 0.568629                              \\ \hline
satellite\_image\_294                         & 0.475034                             & 0.219235                             & 0.221151                             & 0.135448                           & 0.222828                           & 0.584153                            & 0.204638                              \\ \hline
socmob\_541                                   & 0.585336                             & 0.524908                             & 0.510845                             & 0.393824                           & 0.094648                           & 0.53146                             & 0.400179                              \\ \hline
space\_ga\_507                                & 0.353204                             & 0.402592                             & 0.136382                             & 0.21755                            & 0.432612                           & 0.379037                            & 0.243522                              \\ \hline
stock\_223                                    & 0.438271                             & 0.400944                             & 0.336365                             & 0.213751                           & 0.471414                           & 0.580025                            & 0.291354                              \\ \hline
tecator\_505                                  & 0.524734                             & 0.154821                             & 0.717786                             & 0.043018                           & 0.04003                            & 0.653023                            & 0.093625                              \\ \hline
us\_crime\_315                                & 0.306967                             & 0.143614                             & 0.134752                             & 0.018068                           & 0.139102                           & 0.391001                            & 0.086693                              \\ \hline
vinnie\_519                                   & 0.49837                              & 0.414396                             & 0.488711                             & 0.181717                           & 0.330411                           & 0.575668                            & 0.224538                              \\ \hline
visualizing\_galaxy\_690                      & 0.456523                             & 0.285521                             & 0.402726                             & 0.1038                             & 0.064846                           & 0.598021                            & 0.178512                              \\ \hline
wind\_503                                     & 0.399667                             & 0.150179                             & 0.207135                             & 0.113408                           & 0.190706                           & 0.492619                            & 0.208961                              \\ \hline
wine\_quality\_287                            & 0.350762                             & 0.324142                             & 0.156606                             & 0.217912                           & 0.262628                           & 0.42454                             & 0.223012                              \\ \hline
\end{tabular}
\caption{Average empirical stability measured using Jaccard metric.}
\label{jaccard.table}
\end{table*}

\clearpage
\newpage

\begin{table*}[]
\begin{tabular}{|l|r|r|r|r|r|r|r|}
\hline
\multicolumn{1}{|c|}{\textbf{Dataset}}   & \multicolumn{1}{c|}{\textbf{MOSS-H}} & \multicolumn{1}{c|}{\textbf{MOSS-M}} & \multicolumn{1}{c|}{\textbf{MOSS-L}} & \multicolumn{1}{c|}{\textbf{FIRE}} & \multicolumn{1}{c|}{\textbf{GLRM}} & \multicolumn{1}{c|}{\textbf{SIRUS}} & \multicolumn{1}{c|}{\textbf{RuleFit}} \\ \hline
Ailerons\_296                                 & 0.541128                             & 0.225864                             & 0.269797                             & 0.243849                           & 0.230194                           & 0.741071                            & 0.239569                              \\ \hline
ESL\_1027                                     & 0.562471                             & 0.440303                             & 0.523019                             & 0.232895                           & 0.280307                           & 0.683056                            & 0.442317                              \\ \hline
Mercedes\_Benz\_42570 & 0.642431                             & 0.520788                             & 0.453549                             & 0.318953                           & 0.321940                           & 0.695512                            & 0.272624                              \\ \hline
Moneyball\_41021                              & 0.586868                             & 0.347744                             & 0.521799                             & 0.166219                           & 0.255783                           & 0.681284                            & 0.261877                              \\ \hline
abalone\_183                                  & 0.446886                             & 0.368871                             & 0.335289                             & 0.265716                           & 0.36733                            & 0.537703                            & 0.195091                              \\ \hline
autoMpg\_196                                  & 0.568815                             & 0.266127                             & 0.439245                             & 0.22938                            & 0.145534                           & 0.611416                            & 0.169293                              \\ \hline
auto\_price\_195                              & 0.313874                             & 0.168083                             & 0.390849                             & 0.025765                           & 0.054132                           & 0.286174                            & 0.082039                              \\ \hline
bank32nh\_558                                 & 0.615667                             & 0.38057                              & 0.124233                             & 0.235813                           & 0.624428                           & 0.694269                            & 0.418547                              \\ \hline
bodyfat\_560                                  & 0.691549                             & 0.332096                             & 0.69787                              & 0.199768                           & 0.190693                           & 0.744134                            & 0.203302                              \\ \hline
cpu\_small\_227                               & 0.655048                             & 0.439824                             & 0.379241                             & 0.324052                           & 0.469916                           & 0.828172                            & 0.599613                              \\ \hline
elevators\_216                                & 0.485825                             & 0.425721                             & 0.217833                             & 0.332422                           & 0.310246                           & 0.624132                            & 0.412466                              \\ \hline
house\_16H\_574                               & 0.706036                             & 0.337601                             & 0.30926                              & 0.271507                           & 0.612085                           & 0.806731                            & 0.668938                              \\ \hline
houses\_537                                   & 0.760586                             & 0.522144                             & 0.42876                              & 0.580846                           & 0.659909                           & 0.84107                             & 0.691953                              \\ \hline
kin8nm\_189                                   & 0.353694                             & 0.410896                             & 0.190603                             & 0.249942                           & 0.191684                           & 0.387209                            & 0.37927                               \\ \hline
mtp\_405                                      & 0.568483                             & 0.27769                              & 0.159129                             & 0.168782                           & 0.109426                           & 0.657098                            & 0.345234                              \\ \hline
mv\_344                                       & 0.699932                             & 0.640502                             & 0.781384                             & 0.547047                           & 0.651239                           & 0.901553                            & 0.711976                              \\ \hline
no2\_547                                      & 0.305845                             & 0.297938                             & 0.258686                             & 0.173698                           & 0.244973                           & 0.31061                             & 0.294099                              \\ \hline
pol\_201                                      & 0.815246                             & 0.611902                             & 0.639754                             & 0.702624                           & 0.522728                           & 0.878315                            & 0.695217                              \\ \hline
pollen\_529                                   & 0.697261                             & 0.597242                             & 0.347057                             & 0.399166                           & 0.416767                           & 0.697328                            & 0.484955                              \\ \hline
puma32H\_308                                  & 0.665015                             & 0.852149                             & 0.425041                             & 0.387022                           & 0.120502                           & 0.698762                            & 0.715348                              \\ \hline
satellite\_image\_294                         & 0.639354                             & 0.35289                              & 0.347273                             & 0.230399                           & 0.349314                           & 0.731184                            & 0.330443                              \\ \hline
socmob\_541                                   & 0.73334                              & 0.686206                             & 0.664807                             & 0.558244                           & 0.170639                           & 0.687682                            & 0.565676                              \\ \hline
space\_ga\_507                                & 0.514504                             & 0.5696                               & 0.228181                             & 0.352075                           & 0.584631                           & 0.544893                            & 0.386662                              \\ \hline
stock\_223                                    & 0.604512                             & 0.564316                             & 0.496177                             & 0.344409                           & 0.623711                           & 0.726325                            & 0.448297                              \\ \hline
tecator\_505                                  & 0.686373                             & 0.265458                             & 0.830009                             & 0.079659                           & 0.073141                           & 0.781897                            & 0.167463                              \\ \hline
us\_crime\_315                                & 0.448203                             & 0.245703                             & 0.231679                             & 0.034702                           & 0.094120                           & 0.541423                            & 0.155004                              \\ \hline
vinnie\_519                                   & 0.657658                             & 0.581303                             & 0.645009                             & 0.294261                           & 0.488865                           & 0.717735                            & 0.354748                              \\ \hline
visualizing\_galaxy\_690                      & 0.615562                             & 0.43737                              & 0.562095                             & 0.183542                           & 0.106421                           & 0.738236                            & 0.297596                              \\ \hline
wind\_503                                     & 0.564711                             & 0.253707                             & 0.335764                             & 0.192989                           & 0.301108                           & 0.653679                            & 0.339413                              \\ \hline
wine\_quality\_287                            & 0.514151                             & 0.48525                              & 0.260365                             & 0.345662                           & 0.393642                           & 0.581621                            & 0.356523                              \\ \hline
\end{tabular}
\caption{Average empirical stability measured using Ochiai metric}
\label{ochiai.table}
\end{table*}

\clearpage
\newpage

\begin{table*}[h]
\begin{tabular}{|c|c|c|c|c|c|c|c|}
\hline
\textbf{Dataset}    & \textbf{MOSS-H} & \textbf{MOSS-M} & \textbf{MOSS-L} & \textbf{FIRE} & \textbf{GLRM} & \textbf{SIRUS} & \textbf{RuleFit} \\ \hline
Ailerons\_296            & 0.534274        & 0.224652        & 0.262974        & 0.229377      & 0.224151     & 0.732800       & 0.237581         \\ \hline
ESL\_1027                & 0.542947        & 0.424558        & 0.521958        & 0.232340      & 0.276984      & 0.659365       & 0.427253         \\ \hline
Mercedes\_Benz\_42570    & 0.645853        & 0.525220        & 0.437106        & 0.326976      & 0.301945      & 0.697070       & 0.273472         \\ \hline
Moneyball\_41021         & 0.580138        & 0.350952        & 0.525092        & 0.167814      & 0.258333      & 0.673179       & 0.257835         \\ \hline
abalone\_183             & 0.438307        & 0.363968        & 0.345132        & 0.261420      & 0.357088      & 0.530265       & 0.190635         \\ \hline
autoMpg\_196             & 0.591775        & 0.265136        & 0.449240        & 0.221644      & 0.151592      & 0.638246       & 0.173451         \\ \hline
auto\_price\_195         & 0.310888        & 0.167738        & 0.413344        & 0.027336      & 0.052137      & 0.286918       & 0.083685         \\ \hline
bank32nh\_558            & 0.618307        & 0.382116        & 0.127171        & 0.235653      & 0.624348      & 0.697249       & 0.419471         \\ \hline
bodyfat\_560             & 0.697452        & 0.336048        & 0.675731        & 0.195637      & 0.197125      & 0.750574       & 0.205783         \\ \hline
cpu\_small\_227          & 0.651901        & 0.437399        & 0.375417        & 0.320944      & 0.489577      & 0.824803       & 0.599135         \\ \hline
elevators\_216           & 0.504497        & 0.443254        & 0.212563        & 0.321528      & 0.323980      & 0.645238       & 0.432407         \\ \hline
house\_16H\_574          & 0.705332        & 0.338083        & 0.307953        & 0.267354      & 0.605561      & 0.806162       & 0.667867         \\ \hline
houses\_537              & 0.759752        & 0.521477        & 0.427733        & 0.579271      & 0.684691      & 0.840765       & 0.692597         \\ \hline
kin8nm\_189              & 0.354045        & 0.410996        & 0.182477        & 0.245684      & 0.196083      & 0.386821       & 0.377835         \\ \hline
mtp\_405                 & 0.573716        & 0.281380        & 0.162165        & 0.167150      & 0.107181      & 0.663492       & 0.345849         \\ \hline
mv\_344                  & 0.700041        & 0.643118        & 0.776825        & 0.541083      & 0.660033      & 0.901245       & 0.712186         \\ \hline
no2\_547                 & 0.300969        & 0.296931        & 0.266697        & 0.184595      & 0.258225      & 0.306772       & 0.294892         \\ \hline
pol\_201                 & 0.814457        & 0.611534        & 0.650009        & 0.736759      & 0.516296      & 0.878620       & 0.694048         \\ \hline
pollen\_529              & 0.691852        & 0.592593        & 0.353193        & 0.402698      & 0.419920      & 0.692169       & 0.481270         \\ \hline
puma32H\_308             & 0.687810        & 0.878869        & 0.416443        & 0.376662      & 0.123243      & 0.723240       & 0.738950         \\ \hline
satellite\_image\_294    & 0.627831        & 0.346085        & 0.342532        & 0.225709      & 0.358154      & 0.717143       & 0.323611         \\ \hline
socmob\_541              & 0.742848        & 0.687619        & 0.680431        & 0.558221      & 0.173757      & 0.694449       & 0.569620         \\ \hline
space\_ga\_507           & 0.507188        & 0.562345        & 0.245779        & 0.366538      & 0.619276      & 0.535657       & 0.382008         \\ \hline
stock\_223               & 0.617807        & 0.575305        & 0.483883        & 0.344834      & 0.641803      & 0.741020       & 0.457672         \\ \hline
tecator\_505             & 0.680476        & 0.262857        & 0.816622        & 0.079312      & 0.074212      & 0.774074       & 0.164815         \\ \hline
us\_crime\_315           & 0.456817        & 0.251904        & 0.229608        & 0.034443      & 0.140149      & 0.551659       & 0.157081         \\ \hline
vinnie\_519              & 0.661026        & 0.584681        & 0.638624        & 0.288877      & 0.549118      & 0.721775       & 0.357468         \\ \hline
visualizing\_galaxy\_690 & 0.610000        & 0.433968        & 0.570468        & 0.182116      & 0.106315      & 0.731005       & 0.294392         \\ \hline
wind\_503                & 0.580159        & 0.258862        & 0.342267        & 0.190028      & 0.297192      & 0.668519       & 0.348710         \\ \hline
wine\_quality\_287       & 0.502963        & 0.476227        & 0.262332        & 0.347411      & 0.386793      & 0.566679       & 0.350028         \\ \hline
\end{tabular}
\caption{Average empirical stability measured using POG metric.}
\label{pog.table}
\end{table*}
\clearpage
\newpage

\subsection{Additional Discussions}

In the following sections, we provide additional discussions on our MOSS framework.

\subsubsection{\textbf{Generating Candidate Rules}}

Throughout this paper, we use random forests to generate large collections of candidate rules, on which we apply MOSS. Random forests fit decision trees on bootstrapped samples of the original data, where only a subset of features are considered at each split in each tree. The randomness injected from the bootstrap and the feature sub-setting helps create a diverse set of candidate rules.

We note that alternative forms of randomness can be injected into a random forest. For example, in the original random forest paper \cite{breiman2001random}, the author explores randomizing the outputs of each decision tree. Here, we explore adding additional forms of randomness when constructing candidate rules.

We use this new procedure. Given data $X$ and response $y$ we add Gaussian noise to the response to generate $y'$. We then fit a decision tree on $(X,y')$ while randomizing the features considered per split. We repeat until we have a large collection of candidate rules.

Using this procedure, we repeat our experimental setup from \S\ref{experiments.section} and apply MOSS to construct stable rule sets. We show the results in Table \ref{new_randomness.table}.
\vspace{5mm}
\begin{table}[h]
\scalebox{0.9}{
\begin{tabular}{|c|c|c|c|c|c|}
\hline
Method                 & MOSS & FIRE & GLRM & SIRUS & RuleFit \\ \hline
Avg. Accuracy Ranking  & 2.4  & 2.6  & 1.8  & 3.9   & 4.1     \\ \hline
Avg. Stability Ranking & 1.9  & 4.3  & 4.4  & 1.8   & 2.6     \\ \hline
Combined Metric        & 2.2  & 3.5  & 3.1  & 2.9   & 3.5     \\ \hline
\end{tabular}}
\vspace{1mm}
\caption{Experimental results: injecting additional forms of randomness when generating candidate rules.}
\label{new_randomness.table}
\end{table}

These results are consistent with the main results of our paper, and we see here that again MOSS balance predictive accuracy with empirical stability. Only GLRM beats MOSS in terms of accuracy, however, MOSS is significantly more stable; GLRM is the least stable method. Only SIRUS beats MOSS in terms of stability, however, SIRUS is much less accurate. Exploring new methods to generate candidate rules may be an interesting direction for future research.

\subsubsection{\textbf{Classification Tasks}} Our current discussion of MOSS focuses on constructing stable rule sets for regression tasks. However, we can extend MOSS to classification by replacing the ridge-regularized quadratic loss function in the accuracy objective $H_2(z)$ with a ridge-regularized logistic loss function. We could then apply a cutting-plane algorithm similar to the one described in \citep{bertsimas2021sparse}.

\end{document}